\documentclass[11pt]{article}
\usepackage{hyperref}
\usepackage{bbm}
\usepackage{amsmath,amsthm,latexsym,amssymb}
\usepackage{enumerate}
\usepackage{algorithm,enumerate}
\usepackage{algpseudocode}
\usepackage{hyperref}
\allowdisplaybreaks

\usepackage{color}\usepackage{graphicx}

\def\cA{{\mathcal A}}
\def\cB{{\mathcal B}}
\def\cS{{\mathcal S}}
\def\cC{{\mathcal C}}
\def\cG{{\mathcal G}}

\newcommand{\set}[1]{\left\{#1\right\}}

\def\cH{\mathcal{H}}

\def\cP{\mathcal{P}}

\def\ii_(#1,#2){i_{#1}^{#2}}

\def\bx{{\bf x}}

\def\d{\delta}

\def\e{\varepsilon}

\def\z{\zeta}

\def\l{\lambda}

\def\cE{\mathcal{E}}

\def\cS{\mathcal{S}}

\newcommand{\brac}[1]{\left( #1 \right)}

\def\E{{\bf E}}

\renewcommand{\Pr}{\operatorname{\bf Pr}}
\newcommand\bfrac[2]{\left(\frac{#1}{#2}\right)}

\newtheorem{theorem}{Theorem}[section]

\newtheorem{lemma}[theorem]{Lemma}

\newtheorem{observation}[theorem]{Observation}
\newtheorem{remark}[theorem]{Remark}

\newtheorem{notation}[theorem]{Notation}
\date{}

\title{Fast algorithms for solving the Hamilton Cycle problem with high probability}
\author{Michael Anastos \\Freie Universit\"{a}t Berlin \\manastos@zedat.fu-berlin.de }
\date{}

\begin{document}
\maketitle
\begin{abstract}
We study the Hamilton cycle problem with input a random graph $G\sim G(n,p)$ in two settings. In the first one, $G$ is given to us in the form of randomly ordered adjacency lists while in the second one we are given the adjacency matrix of $G$. In each of the settings we give a deterministic algorithm that w.h.p. either it finds a Hamilton cycle or it returns a certificate that such a cycle does not exists, for $p\geq 0$. The running times of our algorithms are w.h.p. $O(n)$ and $O(\frac{n}{p})$ respectively each being best possible in its own setting. \end{abstract}
\maketitle

\section{Introduction}
A Hamilton cycle in a graph is a cycle that visits every vertex exactly ones. The Hamilton cycle problem asks to determine whether a given graph $G$ has a Hamilton cycle or not. It belongs to the list of Karp's 21 $\mathcal{NP}$-complete problems \cite{karp1972} implying that if $\mathcal{P}\neq \mathcal{NP}$ then for every algorithm $\mathcal{A}$ and polynomial $p(\cdot)$ there exists a graph $G$ such that $\mathcal{A}$ cannot solve the Hamilton cycle problem with input  $G$ in $O(p(|V(G)|))$ time.

The problem of finding a fast algorithm for solving the Hamilton cycle problem changes drastically if you require that the algorithm is fast for almost all the graphs instead for all the graphs. Concretely Bollob\'as, Fenner and Frieze \cite{bollobasfennerfrrieze}
gave a deterministic algorithm HAM that takes as input the adjacency matrix of a graph $G$, runs in $O(n^{4+o(1)})$ time and has the property that if $G\sim G(n,p)$ then
$$\lim_{n\to \infty}\Pr( \text{ HAM finds a Hamilton cycle in }G) = \lim_{n\to \infty}\Pr( G\text{ has a Hamilton cycle } ),$$
for all $p=p(n)\geq 0$.
Here by $G(n,p)$ we denote the random binomial graph model i.e. if $G\sim G(n,p)$ then $G$ is a graph on $[n]$  where each edge appears in $G$ independently with probability $p$. The threshold for the existence of Hamilton cycles in $G(n,p)$ is well known. Building upon work of P\'osa \cite{posa1976} and Korshunov \cite{korshunov}, Bollob\'as \cite{bollobas1984}, Koml\'os and Szemer\'edi \cite{komlos1983} and Ajtai, Koml\'os and Szemer\'edi \cite{komlos1983} \cite{ajtai1985} proved that if $G\sim G(n,p)$ where $p=\frac{\log n+\log\log n+c_n}{n}$ then,
\begin{align*}
    \lim_{n\to \infty}\Pr(G \text{ is Hamiltonian})=
    \lim_{n\to \infty}\Pr(\delta(G) \geq 2 )=
    \begin{cases}
0&\text{ if } \lim_{n \to \infty} c_n= -\infty,
\\e^{-e^{-c}} &\text{ if } \lim_{n \to \infty} c_n=c,
\\1&\text{ if } \lim_{n \to \infty} c_n= +\infty.
    \end{cases}
\end{align*}
For a survey on Hamilton cycles in random graphs see \cite{frieze2019}.

The result of Bollobas, Fenner and Frieze was improved, in terms of running time, for some values of $p$, first by Guravich and Shelah \cite{gurevich1987} then by Thomason \cite{thomason89} and lastly by Alon and Krivelevich \cite{alon2020}. All of the 3 corresponding algorithms run in $O(\frac{n}{p})$ time. Observe that in this setting an $O(\frac{n}{p})$ running time is best possible as any algorithm w.h.p.\footnote{We say that a sequence of events $\{\mathcal{E}_n\}_{n\geq 1}$
holds {\em{with high probability}} (w.h.p.\@ in short) if $\lim_{n \to \infty}\Pr(\mathcal{E}_n)=1-o(1)$.} needs to read/query $\Omega(\frac{n}{p})$ entries of the adjacency matrix of $G$ in order to identify $n$ edges, the number of edges in a Hamilton cycle. Ferber, Krivelevich, Sudakov and Vieira show that $(1+o(1))n/p$ queries suffices is general \cite{ferber}.

In another line of research the Hamilton cycle problem for $G\sim G(n,p)$ was consider in the setup where $G$ is given in the form of randomly ordered adjacency lists. That is for every vertex $v\in [n]$ we are given a random permutation of its neighbors via a list $L(v)$. In this setup, Angluin and Valiant \cite{angluin1979} gave a randomize algorithm that finds a Hamilton cycle w.h.p. in $G(n,p)$ for $p \geq \frac{C\log n}{n}$ where C is a sufficiently large constant. Their result was first improved with respect to $p$ first by Shamir \cite{shamir1983}
and then by Bollob\'as, Fenner and Frieze \cite{bollobasfennerfrrieze} whose result is optimal with respect to $p$. Recently it was improved by Nenadov, Steger and Su with respect to the running time \cite{nenadov2020}. We summarize the above results at the table below (taken from \cite{nenadov2020}).
\begin{center}\label{table1}
\begin{tabular}{ l l l l l }
\hline
\textbf{ Authors} &\textbf{ Year} &\textbf{ Time} &\textbf{  p(n) }& \textbf{ Input}  \\ \hline 
 Angluin, Valiant \cite{angluin1979} & `79 & $O(n \log^2n)$ & $p\ge \frac{C_1 \log(n) }{n} $ &  adj. list \\ \hline
 Shamir \cite{shamir1983} & `83 & $ O(n^2)$ & $p\ge \frac{\log n + (3+\epsilon)\log\log n }{n} $ &  adj. list \\ \hline
 Bollobas, Fenner, Frieze \cite{bollobasfennerfrrieze} & `87 &$O( n^{4 + o(1)}) $ & $p\ge 0$ & adj. list or matrix  \\ \hline
 Gurevich, Shelah \cite{gurevich1987}& `87 &$ O(n/p) $ & $p$ const.&  adj. matrix\\ \hline
 Thomason \cite{thomason89} & `89 &$  O(n/p) $ & $p\ge C_2 n^{-1/3} $ &  adj. matrix \\ \hline
 Alon, Krivelevich \cite{alon2020} & `20 &$ O(n/p)$ & $p\ge 70n^{-1/2} $ & adj. matrix \\ \hline
 Nenadov, Steger, Su \cite{nenadov2020} & `21+ & $O(n)$ & $p\geq \frac{C_3\log n }{n}$ & adj. list \\ \hline
\end{tabular}
\end{center}

\subsection{Our results}
In this paper we consider both regimes of how the input is provided and we derive optimal results in both cases.
We say that a sequence of events $\{\mathcal{E}_n\}_{n\geq 1}$
holds {\em{with sufficiently high probability}} (w.s.h.p.\@ in short) if $\lim_{n \to \infty}\Pr(\mathcal{E}_n)=1-o(n^{-7})$.

\begin{theorem}\label{thm:adjlists}
There exists a deterministic algorithm that takes as input a graph $G$ in the form of randomly ordered adjacency lists form and outputs either a Hamilton cycle of $G$ or a certificate that $G$ is not Hamiltonian or FAILURE. If $G\sim G(n,p)$ for some $p=p(n)\geq 0$ then, w.s.h.p. it outputs either a Hamilton cycle of $G$ or a certificate that $G$ is not Hamiltonian in $O(n)$ time.
\end{theorem}
It is worth pointing out that our algorithm is deterministic as opposed to the algorithm given by Nenadov, Steger and Su which is randomized. On the other hand, in \cite{angluin1979}, Angluin and Valiant give a randomized algorithm and show  how it can be derandomized using the randomness found in edges that will not be used by the algorithm in $O(n\log n\log\log n)$ additional time.  

\begin{theorem}\label{thm:adjmatrix}
There exists a deterministic algorithm that takes as input the adjacency matrix of a graph $G$ and outputs either a Hamilton cycle of $G$ or a certificate that $G$ is not Hamiltonian or FAILURE. If $G\sim G(n,p)$ for some $p=p(n)\geq 0$ then,  w.s.h.p. it outputs either a Hamilton cycle of $G$ or a certificate that $G$ is not Hamiltonian by making at most $(1+o(1/\log\log n))n/p$ queries and running in $O(\frac{n}{p})$ time.
\end{theorem}

We use the algorithm that we will present for Theorem \ref{thm:adjmatrix} in the companion paper \cite{anastos2021expected}. There, we present an algorithm that solves the Hamilton cycle problem in $O(\frac{n}{p})$ time when $G\sim G(n,p)$ and $p \geq \frac{C}{n}$. For that, as a first step, we need an Algorithm that with sufficiently high probability solves the Hamilton cycle problem in  $O(\frac{n}{p})$ time. Having probability at least $1-n^{-7}$ suffices.

A $2$-matching of a graph $G$ is a set of edges $M\subseteq E(G)$ such that every vertex in $G$ is incident to at most 2 edges in $M$, as oppose to a matching of size 2.  We say that a 2-matching $M$ saturates a set of vertices $S\subset V(G)$ if every vertex in $S$ is incident to exactly 2 edges in $M$. Thus a necessary condition for $G$ to be Hamiltonian is that for every $S\subset V(G)$ there exists some 2-matching $M_S$ that saturates $S$. The certificates at Theorems \ref{thm:adjlists}, \ref{thm:adjmatrix} will consist of sets of vertices  that violate this necessary condition. 

To prove both Theorems \ref{thm:adjlists}, \ref{thm:adjmatrix} we first build a small collection of vertex disjoint paths that cover the vertex set. We then use rotations-extensions a.k.a.\ P\'osa rotations   to merge these paths into a Hamilton cycle. We present the proofs of Theorems \ref{thm:adjlists} and \ref{thm:adjmatrix} at Sections \ref{sec:adjlists} and \ref{sec:adjmatrix} respectively.

Throughout the paper we omit floors and ceilings whenever these changes do not affect the argument. In addition we assume that $n$ is sufficiently large. Finally in various places we use the Chernoff bounds, stated below, to bound the probability that the binomial random variable $Bin(n,p)$ deviates from its expectation by a multiplicative factor.
\begin{theorem}
For $n\in \mathbb{N}$, $p\in (0,1)$ and $\epsilon>0$,
\begin{equation}\label{eq:chernoffupper}
    \Pr(Bin(n,p)\geq (1+\epsilon)np) \leq \exp\{-\epsilon^2np/(2+\epsilon)\} 
\end{equation}
and
\begin{equation}\label{eq:chernofflower}
    Pr(Bin(n,p)\geq (1-\epsilon)np) \leq \exp\{-\epsilon^2np/2\}.
\end{equation}
\end{theorem}

\section{Hamilton cycles in the randomly ordered adjacency lists model}\label{sec:adjlists}
Let $G\sim G(n,p)$ for some $p=p(n)\geq 0$. In this section we assume that for $v\in [n]$ we are given a randomly ordered adjacency list $L_v$ consisting of the neighbors of $v$ in $G$. We assume that we can query any entry of any list in $O(1)$ time and for $v\in [n]$ and $i\geq 1$ the $i^{th}$ entry of $L_v$, denoted by $L_v(i)$, is ``undefined" if $L_v$ has size strictly less than $i$. 

We let $p^*=\frac{0.999\log n}{n}$ and $q=q(p)=\min\{0.25np,100\log n\}$.

We prove Theorem \ref{thm:adjlists} by implementing a 3-Phase algorithm.
During Phase 1 we  partition $[n]$ into 7 sets $V_1,V_2,...,V_7$. The first 4 will be of equal size while $V_5,V_6$ and $V_7$ will be of sublinear size. In addition we  find a set of vertex disjoint paths $\cP$ from $V_5$ to $V_6$  that cover $V_5\cup V_6 \cup V_7$. Thereafter in Phase 2 we construct perfect matchings in $G[V_i,V_{i+1}]$ for $i=1,2,3$ and in $G[V_1\cup V_5,V_{4}\cup V_6]$ (for disjoint sets $X,Y\subset [n]$ we denote by $G[X,Y]$ the bipartite subgraph of $G$  induced by $X\times Y$). These 4 matchings will result to a 2-factor of $G$. In order to argue that w.s.h.p.\ the resulting 2-factor consists of $O(\log n)$ cycles before generating the last matching we use the randomness found in edges that are not used during this procedure in order to randomly relabel the vertices of $G[V_1\cup V_5,V_{4}\cup V_6]$. In Phase 3 we use extensions-rotations, a.k.a.\ P\'osa rotations, to merge the cycles of the 2-factor into a Hamilton cycle. 

Throughout the 3-Phase algorithm we let $q_i$ be an ordered couple in $[n]\times ([n] \cup \{\text{undefined}\})$ corresponding to the $i^{th}$ query. $q_i=(a_i,b_i)$ if at the $i$th query   an entry of $L_{a_i}$ was queried and the answer was $b_i$. We also let $Q_i=q_1,q_2,...,q_i$ be the sequence corresponding to the first $i$ queries. The algorithm that we consider examines the entries of the lists in order, that is it will never query $L_v(i)$ before querying $L_v(i-1)$ for $v\in [n]$, $i\geq 0$. Thus as the algorithm we consider is deterministic there exists a function $f()$ from sequences with elements in $[n]\times ([n] \cup \{\text{undefined}\})$ to $[n]$ that takes as input $Q_{i-1}$ and returns which list should be queried next,  i.e. $f(Q_{i-1})=a_i$ for $i>0$. For a sequence $X$ with elements in $[n] \times ([n] \cup \{\text{undefined}\})$
we let $D(X)$ be a digraph on $[n]$ where $(a,b)\in E(G(X))$ iff $(a,b)\in X$ for some $b \neq$ undefined and we let $G(X)$ be the simple undirected version of $D(X)$, that is $\{u,v\}\in E(G(X))$ iff $(v,u)$ or $(u,v)$ lies in $E(D(X))$. Finally we let $d_X^+(v),d_X^-(v)$ and $d_X(v)$ be the out-degree, in-degree and the degree of $v$ in $D(X)$, $D(X)$ and $G(X)$ respectively.

Let $\tau=\min\big\{i>0: d_{Q_i}(v)- d_{Q_i}^+(v)>\frac{\log n}{\log\log\log n} \text{ for some } v\in [n] \text{ such that } d_G(v)\geq 3.001q\big\}$. Throughout the algorithm we  ensure that $Q_i$ satisfies the following conditions:
\begin{itemize}
\vspace{-5mm}
    \item[(C1)]$i\leq \min\{10^5n,\tau\}$.
\vspace{-3mm}
\item[(C2)] $L_v(3q +1)$ is not queried for $v\in [n]$.
\vspace{-3mm} 
\item[(C3)] 
$L_v(i-1)$ is queried before $L_v(i)$ is for $v\in [n]$ and $i>0$.
\vspace{-3mm}    \item[(C4)] If $L_v(i)$ is undefined then $L_v(i+1)$ is not queried.
\end{itemize}

Before preceding to the first phase of our algorithm we state a Lemma, which among others, characterizes the distribution of the answer to the next query as long as we have made $O(n)$. Its proof is given at Appendix \ref{app:lists}.

For $i>0$ we let $U_i=\{v\in [n]:(v,\text{ undefined})\in Q_{i}\}$ and $D$ be the set of vertices of degree at most $3.001q$ in $G$.

\begin{lemma}\label{lem:exposerandomness}
If $p\geq p^*$ then w.s.h.p. the following hold:
\begin{itemize}
    \item[(a)] $|D| \leq n^{0.97}$.
    \item[(b)]  For every $v\in[n]$ there exists at most 300 vertices in $D$ from $v$ within distance 10.
    \item[(c)] Let $i< \min\{10^{5}n,\tau \}$ and $v=f(Q_{i})$. If $v\notin D$ then for $u\in [n]\setminus D$,
\begin{align*}
\Pr(q_{i+1}=(v,u)|Q_{i})
\begin{cases}
= (1+o(1)) \cdot \frac{1}{n} \cdot \frac{np-d_{Q_{i}}^+(u)}{np}  &\text{ if } (v,u),(u,v) \notin Q_{i}.
\\  \leq \frac{2000}{np} &\text{ if }  (u,v) \in Q_i
\\ = 0 &\text{ if } (v,u)\in Q_i.
\end{cases}
 \end{align*}
\item[(d)] $\tau > 10^5 n$.
\end{itemize}
\end{lemma}

\subsection{Phase 1: Partitioning the vertex set}
\begin{algorithm}[H]
\caption{}
\begin{algorithmic}[1]
\\Set 
$V_i'=V_i=\{j:(i-1) n/4 <j\leq in/4 \}$ for $i\in[4]$, and $V_5=V_6=V_7=\emptyset$. 
\\Set $I_i=\{n_j: (i-1)n/4\leq j \leq(i-1)n/4+20n^{0.97}\}\subset V_i'$ for $i\in [4]$.
\For{$j\in [4]$ and $v\in V_j$}
\\ \hspace{5mm} Query the entries of $L_v$, one by one, in which of the sets $V_1,...,V_4$ they belong to until either you have identified for each $j'\in [4]$ a set of 4000 entries that belong to $V_{j'}$ or you have queried $min\{q, |L_v|\}$ entries. If the later occurs then remove $v$ from $\cup_{i\in[4]}V_i$ and add it to $V_7$. 
\EndFor
\\ If more than $0.8\cdot 10^5n$ queries have been made in total so far then return $FAILURE0$.
\\  Let $F$ be the graph with edge set $\{\{u,v\}:v\in V_7 \text{ and } u=L_v(j) \text{ for some }j\leq q\}$. 
\For{each component $C$ of $F$ }
\\ \hspace{5mm} Find a set $\cP_C$ of vertex disjoint paths that cover $C \cap V_7$, have no endpoint in $C \cap V_7$ and cover in total at most $3|C \cap V_7|$ many vertices.
\\ \hspace{5mm} If no such set of paths exists then return  FAILURE1 and $S=C\cap V_7$. 
\EndFor
\\ If $p< p^*$ and no FAILURE1 has been returned then return FAILURE2.
\\ Let $\cP=\cup\{\cP_C: \text{ $C$ is a component of $F$}\}$.
\\ For $P\in \cP$ remove $V(P)$ from $\cup_{i\in[4]}V_i$, add one of its endpoints to $V_5$, the other to $V_6$ and all of its interior vertices to $V_7$.
\\Let $d_i= |V_i'\setminus V_i|$ for $i\in[4]$ and  $d_{max}=\max_{i\in[4]} d_i$. 
\\ If $d_{max}> 6n^{0.97}$ then return FAILURE3. 
\\ Else remove from $V_i$ a set of vertices of size $d_{max}-d_i-\mathbb{I}(|V_i'|>0.25n)$ in $V_i\cap I_i$, for $i\in[4]$ and let $U$ be the set of vertices that have just been removed. 
\For{$v\in U$}
\\ \hspace{5mm} Greedily match $v$ to vertices   $w(v,2)\in V_{2}, w(v,3) \in V_{3}$ and then match the vertices $w(v,2)$ and $w(v,3)$  with  vertices $w(v,1) \in V_1$ and $w(v,4) \in V_4$ respectively such that the $|U|$ resultant paths of length 4 are vertex disjoint (the edges used to match vertices in this step are originated from queries that took place at lines 3 to 5).
\\ \hspace{5mm} If  no such vertices exists return FAILURE4.
\EndFor
\\ For each path $P$ created remove $V(P)$ from $\cup_{i\in[4]}V_i$, add one of its endpoints to $V_5$, the other to $V_6$ and all of its interior vertices to $V_7$.
\end{algorithmic}
\end{algorithm} 

\begin{observation}\label{rem:failure1}
If Algorithm1 returns FAILURE1 and $S\subset V(G)$  then $G$ is not Hamiltonian and $S$ is a certificate of this fact. Indeed if such a set $S$ is returned then there does not exist a set of vertex disjoint paths with no endpoints in $S$ that covers $S$. On the other hand any Hamilton cycle of $G$ induces such a set of paths on $S\cup N(S)$. 
\end{observation}
 
\begin{lemma}\label{lem:failure2}
W.s.h.p. Algorithm 1 does not return FAILURE0 or FAILURE2.
\end{lemma}
\begin{proof}
For $i \leq 0.8 \cdot 10^{5}n$ let $Y_i=1$ if $(q_i)_2$ belongs to $V_j$ and for at most 3999 elements in $\{(q_\ell)_2:\ell<i\}$ belong to $V_j$. Lemma \ref{lem:exposerandomness} implies that $\Pr(Y_i=1)\geq 0.24$ independently of $Y_1,Y_2,...,Y_{i-1}$. In addition in the event that Algorithm 1 returns FAILURE0 among the variables $Y_1,Y_2,...,Y_{0.8\cdot 10^5 n}$ at most $16000n$ equal to 1. Thus \eqref{eq:chernofflower} implies,
\begin{align*}
    \Pr(FAILURE0)\leq \Pr(Bin(0.8\cdot 10^5n,0.24)\leq 16000n)\leq \exp\{-\Omega(n)\}=
    o(n^{-7}).
\end{align*}

Algorithm 1 returns FAILURE2 only if $p< p^*$ and $G$ has minimum degree at least 1. By sequentially query $L_j(1)$ for $j\in [n]$ if $j$ does not appear in $\{L_i(1):1\leq i<j\}$ we have,
\begin{align*}
    Pr&( \text{ Algorithm 1 returns FAILURE2 }) \leq (1-(1-p^*)^n)^{0.5n} \leq \bigg(1-e^{-\frac{np^*}{1-p^*}}\bigg)^{0.5n}
=o\big( n^{-7} \big).
\end{align*}
\end{proof} 

\begin{notation}
We call a vertex $v\in [n]$ \emph{ burned} if it enters the set $V_7$ at line 4 of Algorithm 1. We denote the set of burned vertices by $B$.
Thereafter for $v\in [n]$ we let $index_1(v)$ be the first entry of $L_v$ that has not been queried yet and $G_1$ be the graph on $[n]$ whose edges are the ones revealed by Algorithm 1.
\end{notation}
The following Lemma describes typical properties of $B$ and $G_1$.Its proof is located at Appendix \ref{app:lists}.
\begin{lemma}\label{lem:phase1prep}
W.s.h.p. if $p^* \leq p \leq \frac{1000\log n}{n} $ then the following hold:
\begin{itemize}
    \item[(a)] $|B\setminus D|\leq  n^{0.97}$.
    \item[(b)] For every $v\in[n]$ there exists at most 300  vertices in $B\setminus D$ within distance 10 from $v$ in $G_1$.
    \item[(c)] No subgraph $F$ of $G_1$ has at most $10^5$ vertices and spans more than $|V(F)|+8$ edges.  
    \item[(d)] For every $v\in V$ the number of vertices in $\cup_{i\in[4]}I_i$ within distance  at most 10 from $v$ is at most 300.
\end{itemize}
In the case $p\geq \frac{1000\log n}{n}$ we have that w.s.h.p. $B=D=\emptyset$.
\end{lemma}

The following Lemma describes the performance and the correctness of Algorithm 1.

\begin{lemma}\label{lem:alg1}
W.s.h.p. if $p>p^*$ then either Algorithm 1 returns FAILURE1 or the following hold:
\begin{itemize}
    \item[(a)] For  every vertex $v\in[n]$ and $i\in [4]$ at most $1000$ neighbors of $v$ belong to $V_i'\setminus V_i$.
    \item[(b)]  Algorithm 1 does not return FAILURE3.
    \item[(c)]  Algorithm 1 does not return FAILURE4.
    \item[(d)] Algorithm 1 runs in $O(n)$ time.
\end{itemize}
\end{lemma}
\begin{proof}
For (a),(b),(c) we first examine the case $p^*\leq p\leq \frac{1000\log n}{n}$.

(a) Let $v\in [n]$ and $R_u$ be the set of vertices in $B\cup I$ within distance 3 from $v$. As $B\cup I\subset D\cup (B\setminus D) \cup I)$ Lemmas \ref{lem:randomness} and \ref{lem:phase1prep} imply that $|R_u| \leq 900$. Thereafter, as for each  neighbor $u$ of $v$ in $V_i'\setminus V_i$ there exists a path of length at most 2 from $u$ to $R_v$, Lemma \ref{lem:phase1prep} implies that w.s.h.p $v$ has at most $900+9<1000$ such neighbors.

(b) Lemmas \ref{lem:exposerandomness}, \ref{lem:phase1prep} imply that w.s.h.p. $|B| \leq 2n^{0.97}$. Let $Y=V_5\cup V_6\cup V_7$. At line 7 of Algorithm 1 every vertex in  $Y$ belongs to $B$ while from line 7 to line 16 the size of $Y$ may increase by a factor of 3. Hence at line 16 w.s.h.p. $d_i\leq 3|B|\leq 6n^{0.97}$ and Algorithm 1 does not return FAILURE3

(c) Let $v\in U$. Part (a)  implies that w.s.h.p. there exists at least $3000$ neighbors in  the set $V_2$ ($V_3$ respectively) each having at least 3000 neighbors in $V_1$ (resp. $V_4$). In addition Lemma \ref{lem:phase1prep} implies that w.s.h.p. there exists at most 300 vertices in $U$ within distance 10 from $v$. Thus when choosing $w(v,2)$ (resp. $w(v,3)$) there exist at least 3000-300 valid candidates. 
Reiterating the argument gives that when choosing $w(v,1)$ (resp. $w(v,4)$) there exist at least 2700 valid candidates too.

Now if $p\geq \frac{1000\log n}{n}$ Lemma \ref{lem:phase1prep} states that w.s.h.p. $B=D=\emptyset$. Then w.s.h.p. $d_{max}=0$, $|V_j'\setminus V_j|\leq 1$ for $j\in [4]$, $|U|\leq 4$ and (a),(b),(c) hold.

(d) The runtime of the Algorithm 1 is dominated by the number of queries at lines 3 to 5  and the complexity of lines 8 to 11. For each component $C$ of $F$ the set $C\cap V_7$ consists of burned vertices, hence by 
Lemma \ref{lem:phase1prep}, w.s.h.p. $F$ has at most $2n^{0.97}$ components each spanning at most 300 vertices of degree at most $q\leq 100\log n$ in $G_1$. To cover $C\cap V_7$ with suitable paths we can first cover the set of vertices 
of degree at most 600 in $C\cap V_7$ (hence in $G$), by examining all of the at most $2^{300\cdot 600}$ related options, and then greedily cover the  rest of the vertices in $C\cap V_7$ by paths of length 2. Thus w.s.h.p. lines 8 to 11 take O(n) time to be executed and Algorithm 1 runs in $O(n)$ time.
\end{proof}

\subsection{Phase 2: Finding perfect matchings between the partitions}
Let $U_1=V_1 \cup V_5$ and $U_4=V_4\cup V_6$. We apply Algorithm 2 to each of the pairs $(V_1,V_2),(V_2,V_3),$ $(V_3,V_4)$ and $(U_1,U_4)$ to find a perfect matching in the underlying bipartite graph. Before applying the algorithm we order the vertices in $V_i\subset [n]$ in an increasing order for $i\in [4]$. For $i\in\{1,4\}$ we also order the vertices in $U_i$ as follows. For $v\in U_1$ let $f(v)=(v_1,v_2,...,v_{10})$ be the first   10 neighbors of $v$ in $V_{3}$ as they appear in $L_v$. These 10-tuples can be deduced by the queries made during Algorithm 1. Order the vertices in $U_1$ such that $u$ precedes $v$ if $f(u)$ is lexicographical smaller than $f(v)$. Since the 10-tuples $f(v):v\in U_1$ are independent and distributed uniformly at random over 10-tuples with elements in $V_{3}$ we have that in the event $f(u)\neq f(v)$ for each pair of distinct vertices $u,v\in U_1$ this procedure produces a uniform random ordering of $U_1$.
This occurs with probability at least 
$$ \binom{n}{2} \binom{(1+o(1))n/4}{10}^{-1}=o(n^{-7}).$$
Similarly reorder the vertices in $U_4$. As the edge sets used to generate these 2 random orderings are independent, the 2 new orderings of the vertices of $U_1,U_4$ are also independent.

Algorithm 2 takes as input sets $X=\{x_1,x_2,...,x_r\}$, $Y=\{y_1.y_2,...,y_r\}$ of size $r$ and the adjacent lists $L_v, v\in X\cup Y$. For each $v\in X$ (resp. $y\in Y)$ from Algorithm 1 we are given a list $L_v^1$ of the first (defined by $L_v$) 3000 neighbors in $Y$ (resp. in $X$) and the first entry of $L_v$ that has not been queried by Algorithm 1 which we denote by $index_v$. We denote by $H_1$ the graph induced by $L_v^1,v\in X\cup Y$
and for $S\subset X\cup Y$ by $H_1(S)$ the subgraph of $H_1$ induced by $S$.

Algorithm 2 starts by greedily constructing a matching $M$  which it then greedily augment using alternating paths of length 3 to a matching $M''$.
Thereafter, for every unsaturated by $M''$ vertex $w\in X,Y$ it attempts to generates a tree on 2001 vertices rooted at $w$ with exactly 1000 leaves all of which are found at depth 2 such that each root to leaf path is an $M''$-alternating path. The corresponding  trees will be vertex disjoint. The for-loop at line 11 takes as input the trees generated at line 8 and runs for $\log_{10^4} n$ iterations. At the $i^{th}$ iteration it is given the current matching $M_{i-1}$ the set of unsaturated vertices $X_{i-1}\subset X$, $Y_{i-1} \subset Y$, which will be of size
$O(10^{-4i}n)$ and a sets of vertex disjoint trees $T_{i-1}$. $|T_{i-1}|=|X_{i-1}\cup Y_{i-1}|$. Each tree $T$ in $T_{i-1}$ is rooted at some vertex $r\in X_{i-1}\cup Y_{i-1}$
has $O(10^{3i})$ leaves and every root to leaf path of $T$ is $M_{i-1}$-alternating. The algorithm first attempts to augment $M_{i-1}$ by matching leaves and using the underlying $M_{i-1}$-alternating paths. Then it continuous by further growing the trees whose roots are still unsaturated by the current matching.

Throughout the execution of lines 10 to 15 of Algorithm 2 the set $Q$ consists of all the vertices $v$ for which $L_v(index_v)$ may has been queried. The subroutines \textbf{Augment} and \textbf{Extend} are given in a later point of this Subsection. 

\begin{algorithm}[H]
\caption{}
\begin{algorithmic}[1]
\\Let $M'=\emptyset$.
\\While there exists $e\in E(H_1)$ not incident to $M'$, add $e$ to $M'$.
\\Let $X'\subseteq X,Y'\subseteq Y$ be the sets of $M'$-unsaturated vertices.
\While{ there exists $x\in X'$, $y\in Y'$, $v\in L_x^1$ and $u\in L_y^1$ such that $uv\in M'$ }
\\ \hspace{5mm} Remove $x$ from $X'$, $y$ from $Y'$ and set $M'=M'\triangle\{xv,vu,uy\}$. 
\EndWhile 
\\Set $M''=M'$ and let $X''\subset X,Y''\subset Y$ be the sets of $M''$-unsaturated vertices.    
\\$M_{0},F_X^{0},F_Y^{0},FAILURE5=\textbf{GenerateTrees}(M'',X'',Y'',\{L_v^1\}_{v\in X''\cup Y''})$. 
\\If FAILURE5=1 or $|X''|>10^{-6}n$ then return FAILURE.
\\Let $Q=\emptyset$.
\For{ $i=1$ to $\log_{10^4}n$ }
\\ \hspace{5mm} $M_i,F_X^{i-1},F_Y^{i-1},Q=\textbf{Augment}(M_{i-1},F_X^{i-1},F_Y^{i-1},Q,\{L_v\}_{v\in X\cup Y})$.
\\ \hspace{5mm} $ M_i,F_X^i,F_Y^i,FAILURE6, FAILURE7=\textbf{Extend}(M_{i},F_X^{i-1},F_Y^{i-1},Q,\{L_v\}_{v\in X\cup Y},i)$.
\\ \hspace{5mm} If  $|F_X^{i}|>|10^{-4-4i}n|$ or FAILURE6=1 or FAILURE7=1 then return FAILURE.
\EndFor   
\\Output $M=M_{\log_{10^4} n}$.
\end{algorithmic}
\end{algorithm} 

For $S\subset X\cup Y$ let $H_1^L(S)$ be the subgraph of $H_1$ with edge set $\{\{v,u\}:v\in S, u\in L_v^1\}$.
\begin{algorithm}[H]
\caption{${GenerateTrees}(M'',X'',Y'',\{L_v^1\}_{v\in X''\cup Y''})$}
\begin{algorithmic}[1]
\\ Let $F_X=F_Y=\emptyset$, FAILURE5=0 and set $\tilde{X}=X''$ and $\tilde{Y}=Y''$.
\For{ $W\in \{\tilde{X},\tilde{Y}\}$ }
\While{ $E(H_1^L(W))\neq \emptyset$ and there exists $w\in W$ that is incident to a set $S_w$ of 1000 vertices of degree 1 in $H_1^L(W)$ }
\\ \hspace{10mm} Remove $w$ from $W$ and add the tree $T_w$  to $F_W$. $T_w$ is a tree on $2001$ vertices rooted at $w$ of height 2. In $T_w$ the neighbors of $w$ are the vertices in $S_w$ while every vertex in $S_w$ in addition to $w$ is connected to a unique leaf of $T_w$ via an edge in $M''$.
\EndWhile
\EndFor
\\ If $\tilde{X}\neq \emptyset$ or $\tilde{Y}\neq \emptyset$ then set FAILURE5=1
\\ Return $M'',F_X,F_Y,FAILURE5$.
\end{algorithmic}
\end{algorithm} 

\begin{observation}
For $w\in X''\cup Y''$ all of the neighbors of $w$ in $S_w$ are saturated by $M''$ since otherwise the corresponding edge would had been used to augment $M''$ earlier at line 2 of Algorithm 2. 
\end{observation}
The proofs of parts (a)-(c) of the following Lemma are found at Appendix \ref{app:lists}. \begin{lemma}\label{lem:phase2prep}
W.s.h.p. the following hold:
\begin{itemize}
    \item[(a)] Every pair of sets $X'\subseteq X, Y'\subseteq Y$ of size at least $0.002n$ spans an edge in $H_1$.
    \item[(b)] For every $S\subseteq X$ or  $S\subseteq Y$, $|S| \geq 10^{-6}n$ implies that the neighborhood of $S$ in $H_1$ has size at least  $ 0.002n$.
    \item[(c)] For every $S\subset X$ (resp. $S\subset Y$) of size at most $10^{-6}n$ the set $\{u: u\in L_v^1 \text{ for some }v \in S\}$ has size at least $2000|S|$.
    \item[(d)]  $|X''|=|Y''|\leq 10^{-6}n$.
    \item[(e)] FAILURE5=0
\end{itemize}
\end{lemma}

\begin{proof}
(d) Let $N_{X''}=\{v\in L_x^1:x\in X''\}$ and $N_{Y''}=\{u\in L_y^1:y\in Y''\}$. Due to the execution of Algorithm 2 there is no edge from $N_{X''}$ to $N_{Y''}$ in $H_1$ and therefore, by part (a), w.s.h.p. $|N_{X''}|\leq 0.002n$ or $|N_{Y''}| \leq 0.002n$. W.l.o.g assume that the former applies. Then, part (b) implies that w.s.h.p.  $|Y''|=|X''|\leq 10^{-6}n$. 

(e) For $W\in \{\tilde{X},\tilde{Y}\}$ and $S\subset W$ each vertex in $S$ has exactly 3000 neighbors in $H_1^L(S)$. On the other hand part (d) implies that w.s.h.p. there are at least $2000|S|$ distinct vertices incident to $|S|$ in $H_1^L(S)$. Thus there exists at least 1 vertex $s\in S$ that is incident to at least 1000 vertices of degree 1 in $H_1^L(S)$. Hence at each iteration of the while-loop if $W\neq \emptyset$ then a new tree is constructed.
\end{proof}
\begin{algorithm}[H]
\caption{${Augment}(M_{i-1},F_X^{i-1},F_Y^{i-1},Q,\{L_v\}_{v\in X\cup Y})$}
\begin{algorithmic}[1]
\For{ $T\in F_X^{i-1}$} let $L(T)$ be the leaves of $T$.
\For{ $u\in L(T)$}
\\ \hspace{10mm} Reveal the first 1000 entries of $L_{u}$ starting with $L_{v}(index_{u})$ and add $u$ to $Q$. 
\If{$w\in W$ and any entry revealed is a leaf $w$ of a tree $T'\in F_Y^{i-1}$}
\\ \hspace{15mm} Use the path from the root of $T$ to $u$ to $w$ to the root of $T'$ to augment $M_{i-1}$. 
\\ \hspace{15mm} Remove $T$ from $F_X^{i-1}$ and $T'$ from $F_Y^{i-1}$.
\EndIf
\EndFor
\EndFor
\\Return $M_{i-1},F_X^{i-1},F_Y^{i-1},Q$.
\end{algorithmic}
\end{algorithm} 

\begin{algorithm}[H]
\caption{${Extend}(M_{i},F_X^{i-1},F_Y^{i-1},Q,\{L_v\}_{v\in X\cup Y},i)$}
\begin{algorithmic}[1]
\\Let $F_X^{i}=\emptyset, F_Y^{i}=\emptyset, FAILURE6^i=FAILURE7^i=0$.   
\For{$(W,Z) \in \{(X,Y),(Y,X)\}$ and $T\in F_W^{i-1}$}   
\\\hspace{5mm} 
Let $P(T)$ be the set of parents of leaves of $T$. Set $A_T=\emptyset$.
\For{$p\in P(T)$}
\For{a descendent $u$ of $p$ in $T$}
\\ \hspace{10mm} Let $j=index_u+1000$ and $A_u=\emptyset$.
\While{($|A_T|+|A_u|<  10^{3+3i}$ or $|A_u|<1000$)  and $j<index_u+q$ and $L_u(j-1)\neq undefined$}
\\ \hspace{20mm} Query $L_u(j)$. If (i) $L_u(j)$ belongs to $Z$, (ii) $L_u(j)$ does not belong to some tree in $F_X^{i-1}\cup F_X^{i} \cup F_Y^{i-1} \cup F_Y^{i}$, (iii) $L_u(j)\notin e$ for some $e\in A_T$ and (iv) the neighbor of $L_u(j)$ induced by $M_{i}$ does not lie in $Q$ then add  $(u, L_u(j))$ to $A_u$.
\EndWhile
\\\hspace{15mm} If $|A_u|\geq 1000$ then set $A_T=A_T\cup A_u$.
\EndFor
\EndFor
\If{ $|A_T|<10^{3+3i}$}{ set $FAILURE6^i=1$.}
\Else{ extend $T$ by adding to it first the edges in $A$ and then the edges in $M_{i}$ incident to $A_T$. Thereafter delete from $T$ all leaves at depth $2i$. Finally add $T$ to $F_X^i$}
\EndIf
\EndFor
\If{ more than $2^{-i-1}n+q|D|$ queries were made in total }{ set $FAILURE7^i=1$.} 
\EndIf
\\Return $M_{i},F_X^i,F_Y^i,FAILURE6^i, FAILURE7^i$.
\end{algorithmic}
\end{algorithm} 

\begin{observation}\label{obs:3}
Let $T\in F_X^i\cup F_Y^i$ and assume Algorithm 2 did not return FAILURE at the $(i-1){th}$ iteration of the for-loop at line 11 or $i=1$. Then $T$ has size at most $10^{4+3i}$,  height $2i+2$ and $10^{3+3i}$ to $10^{3+3i}+1000$  leaves all at depth $2i+2$. In addition for each leaf there exists at least 999 other leaves at its 4th neighborhood. Thus Lemma \ref{lem:exposerandomness} implies that w.s.h.p. $T$ has at least $0.7\cdot10^{3+3i}$ leaves not in $D$. In addition the trees in $\bigcup_{i\geq 0} (F_X^i\cup F_Y^i)$ span at most $2\cdot2001\cdot10^{-6}n+2\sum_{i\geq 1}(10^{-4-4i} n) ( 10^{3+3i}+1000) < 0.1n$ many vertices. As $Q$ is a subset of those vertices, thought the algorithm $|Q| \leq 0.1n$ and at any point of the algorithm  there exists at most $0.1n$ vertices violating the condition at line 6 of Extend. 
\end{observation}
\begin{lemma}\label{lem:randomness}
W.s.h.p. no FAILURE is return at the $i^{th}$ execution of line 14 at Algorithm 2 for $1\leq i\leq \log_{10^4}n$.
\end{lemma}    
    
\begin{proof}
By Lemma \ref{lem:phase2prep} w.s.h.p. $|F^X_{0}|<10^{-4}n$. For $1\leq i \leq \log_{10^4} n$ and $\cE_i$ be the event that  Algorithm 2 did not return FAILURE at the $j^{th}$ iteration of the for-loop at line 11 for $j\leq i-1$ and $|D|\leq n^{0.97}$. For part (a)  it suffices to show that for $1\leq i\leq \log_{10^4} n$, conditioned on $\cE_i$,  with probability $1-o(n^{-7}\log^{-1}n)$  we have that 
$|F_X^i|\leq 10^{-4-4i}n$ and $FAILURE6^i=FAILURE7^i=0$ at the end of the $i^{th}$ execution of \textbf{Extend}. 

Let $1\leq i\leq \log_{10^4} n$. In the event $\cE_i$ we have that $|F^X_{i-1}|, |F^Y_{i-1}|\leq 10^{-4-4(i-1)}n$  and by Observation \ref{obs:3} each tree in $F^X_{i-1}\cup F^Y_{i-1}$
spans a set of at least  $0.7\cdot 10^{3+3(i-1)}$ leaves not in $D$. During the execution of \textbf{Augment}, while $|F^X_{i-1}|>10^{-4-4i+j}n$, $j\in\{0,1,2,3\}$ the trees in $F^Y_{i-1}$ spans at least  $n_{i,j}=0.7\cdot 10^{3+3(i-1)}\cdot 10^{-4-4i+j}n=0.7\cdot 10^{-4-i+j}$ leaves not in $D$. At the same time a tree  $T\in F^X_{i-1}$ has at least $0.7\cdot 10^{3+3(i-1)}$ many leaves not in $D$ each having $10^3$ entries of its list queried. Lemma \ref{lem:exposerandomness} implies that with probability at least $\frac{(1+o(1))n_{i,j}}{n}$ each of the at least $0.7\cdot 10^{3+3(i-1)+3}$ queries may results to an augmentation of $M_{i-1}$. Thus the root of $T$ will be saturated by $M_i$ with probability at least 
$$1-\bigg(1- \frac{(1+o(1))n_{i,j}}{n} \bigg)^{0.7\cdot 10^{3+3i}}\geq 1-e^{0.49\cdot 10^{2i+j-1}}\geq 1-e^{4.9\cdot 10^{j}}.$$
Hence, \eqref{eq:chernoffupper} gives,
\begin{align*}
\Pr\big(|F^X_i|> 10^{-4-4i}n|\cE_i\big)\leq \sum_{j=0}^3\Pr\big(Bin\big( 10^{-4-4(i-1)+j}n,e^{-4.9\cdot 10^{j}}\big)>0.25 \cdot 10^{-4-4i}n\big)=o(n^{-8}).    
\end{align*}
Now let   $(W,Z)\in\{(X,Y),(Y,X)\}$, $T\in F_W^{i-1}$ and $v_1,v_2,....,v_r$ be the leaves of $T$ in the order that they are examined by \textbf{Extend}. Call a leaf $v_j$ \emph{good} if $|A_{v_j}|\geq 10^4$ and $j<r$ or $v_j$ is not examined by \textbf{Extend}. Thus for $j<r$ if $v_j \notin D$ as $|Z\setminus Q|\geq 0.14n$, Lemma \ref{lem:exposerandomness} implies that  
\begin{align*}
\Pr(v_j \text{ is not good }|\cE_i)&\leq \Pr(Bin(q-1000,0.14)<10^4) )\leq \binom{q}{10^4}0.86^{q-2\cdot 10^4}
\\&\leq (100\log n)^{2\cdot 10^4} 0.86^{0.24\log n}=o(n^{-0.03}). 
\end{align*}
If $T$ has at least $0.1\cdot 10^{3+3(i-1)}$ many good leaves then $|A_T| \geq \min\{10^{3+3i}, 10^{3i-1}\cdot 10^4\}\geq 10^{3+3i}$. By Observation \ref{obs:3} $T$ has at least $0.7\cdot 10^{3+3(i-1)}$ many leaves in $D$. Therefore, 
$$\Pr(|A_T|< 10^{3+3i}|\cE_i) 
\leq \binom{0.7\cdot 10^{3+3(i-1)}}{0.6\cdot 10^{3+3(i-1)}} n^{-0.03 \cdot 0.6\cdot 10^{3+3(i-1) }} =O(n^{-18})$$
and hence $\Pr(FAILURE6^i=1|\cE_i)=O(n^{-17})$.

Now call a leaf $v_j$ \emph{unhelpful} if $|A_{v_j}|<1000$ and  $v_j$ is examined by \textbf{Extend}. Similarly to before,
$\Pr(v_j \text{ is unhelpful}|\cE_i)=o(n^{-0.03})$.
For $i\leq \log\log n$, in the event $\cE_i$ if $FAILURE7^i=1$ then either more than $n^{0.98}$ leaves not in $D$ are unhelpful or among the first $2^{-1-i}n$ queries made by \textbf{Extend}, of entries of lists of vertices not in $D$, at most $(10^{3+3i}+1000)\cdot2 \cdot |F_X^{i-1}|+1000n^{0.98} \leq 3 \cdot 10^{-1-i}n$ satisfy the conditions of line 8. Thus, by \eqref{eq:chernofflower},
\begin{align*}
    \Pr(Failure 7^i=1|\cE_i) &\leq
    \binom{n}{n^{0.98}}n^{-0.03\cdot n^{0.98}}+
    \Pr(Bin(2^{-i-1}n-q|D|,0.14)< 3\cdot 10^{-1-i})
=o(n^{-8}).
\end{align*}

On the other hand, for $\log\log n<i \leq \log_{10^4} n$, in the event $\cE_i$ if $FAILURE7^i=1$ then there exists a tree such that at least $\log^{-2}n$ fraction of its leaves not in $D$ are unhelpful. Thus, with $n_l=0.7\cdot 10^{3+3(i-1)}\geq 10^{3\log\log n} \geq \log^5 n$ and $n_u=10^{3+3(i-1)} +1000$
by \eqref{eq:chernofflower},
\begin{align*}
    \Pr(Failure 7^i=1|\cE_i) &\leq
    2n \sum_{j=n_l}^{n_u} \binom{j}{\frac{j}{\log^2n}}n^{-0.03\cdot \frac{j}{\log^2n}}
\leq     2n \sum_{j=n_l}^{n_u} \big(en^{-0.03}\log^2 n\bigg)^{-\frac{j}{\log^2n}}
\\ &\leq 2n^2 \big(en^{-0.03}\log^2 n\bigg)^{-\frac{n_l}{\log^2n}}
\leq 2n^2 \big(en^{-0.03}\log^2 n\bigg)^{-\log^3n}=o(n^{-8}).
\end{align*}
\end{proof}    
Lemmas \ref{lem:phase2prep} and \ref{lem:randomness} imply that w.s.h.p. Algorithm 2 does not return FAILURE, makes at most $\sum_{i\geq 1}2^{-i-1}n\leq n$ queries and runs in $O(n)$ time in total. In addition it outputs $M=M_{\log_{10}n}$ and the number of $M$-unsaturated vertices is at most
$$2|F_X^{\log_{10^4} n}| \leq 2\cdot 10^{-4-4\log_{10^3}n}<1.$$
Hence $M$ is a perfect matching.

\subsection{Phase 3: Merging the cycles}
Let $F$ be the 2-factor constructed by Algorithm 2. As $F$ induces a random permutation on the vertices of $V_1$ and every cycle of $F$ has such a vertex it is easy to show that w.s.h.p. $F$ consists of $\ell \leq 100\log n$ cycles spanning $[n]$. In addition, due construction every cycle spans an edge from $U_1$ to $U_4$, hence an edge whose endpoints are of degree at least 16000 in $G_1$. 
We start by removing from every cycle such an edge. Let $\cP=\{P_1,P_2,\ldots,P_{\ell+1}\}$ be set of vertex disjoint paths created and for $P_i\in cP$ be $v_{(i,1)}$ and $v_{(i,2)}$ be its endpoints. The cycle $H_0=v_{(1,1)},P_1,v_{(1,2)}v_{(2,1)},P_2,v_{(2,2)}v_{(3,1)}P_3,\ldots,$ $v_{(\ell+1,1)}P_{\ell+1} $ $ v_{(\ell+1,2)}v_{(1,1)}, v_{(1,1)}$ is a Hamilton cycle in the graph $\Gamma_0=(V,E(H)\cup R)$ where $R=\set{\set{v_{(i,2)},v_{(i+1,1)}}:i\in[\ell]} \cup \set{\set{v_{(\ell+1,2)},v_{(1,1)}}} $. We let $V_R=\underset{e\in R}{ \cup}e$. The following Lemma, proven in Appendix \ref{app:lists}, implies that for $v\in R$ there exists sets $S(v),T(v),E(v)$ of size $10^5,3\cdot 10^5$ and $2\cdot 10^5$ respectively such that $S(v)\subset N_{G_1}(v)$, $T(v)$ consists of the $H_0$-neighborhood of $S(v)$  and $E(v)$  consists of the edges in $E(H_0)$ spanned by $T(v)$. Its proof is given in Appendix \ref{app:lists}. 

\begin{lemma}\label{lem:endpoints}
W.s.h.p. there exists sets $S(v)\subset N_{G_1}(v)$ for $v\in R$ such that $|S(v)|\geq 10^4$ for $v\in R$ and no vertices in $V_R\cup \big(\cup_{v\in R} S(v) \big)$ are within distance 10 in $H_0$.
\end{lemma} 
Let $T_R=\cup_{v\in V_R} \big(S(v) \cup T(v) \big)$.

Starting with $H_0$, we find a Hamilton cycle in $G$ by removing from $H_0$ the edges in $R$  one by one. We do this in at most $\ell$ rounds of an extension-rotation procedure. Fix $i\geq 1$ and suppose that after $i-1$ rounds, we have a Hamilton cycle $H_{i-1}$ in $G_{i-1}\subset G\cup R$ such that $|E(H_{i-1})\cap R| \leq \ell-(i-1)$ and a set $Q$ of vertices whose lists have been queried during Phase 3. $Q$ will have size at most $(i-1)n^{0.56}$. In addition suppose that if $e\in R$ belongs to $H_{i-1}$ then for $v\in e$ (i) the edges in $E(v)$ belong to $H_{i-1}$ and (ii) $Q\cap T(v)=\emptyset$. 

We start round $i$ by deleting an edge $e\in R$ from $H_{i-1}$  to create a Hamilton path $P=v_1,v_2,....,v_n$. Given a path $P'=x_1,x_2,\ldots,x_n$ and an edge $\set{x_s,x_n}$ where $1\leq s<n$, the path $x_1,\ldots,x_{s-1},x_s,$ $x_{n},x_{n-1},\ldots,x_{s+1}$ is said to be obtained from $P$ by a P\'osa rotation with $x_1$ as the fixed endpoint. We call the edge $\set{x_s,x_n}$,the edge $\set{x_s,x_{s+1}}$, the vertex $x_j$ and the vertex $x_{j+1}$ the inserting edge, the deleted edge, the pivot vertex and the new endpoint respectively.

We then split $P$ into 3 vertex-disjoint paths  $P_1,P_2,P_3$ each of length at least $n/3-2$. 
We let $j_1,j_n$ be such that  such that $|S(v_1) \cap V(P_{j_1})|$ and $|S(v_n) \cap V(P_{j_n})|$ respectively are maximized. If $j_1\neq j_n$ then let $j_3=[n]\setminus \{j_1,j_n\}$, $J_1= (V(P_{j_1})\cup V(P_{j_3}))\setminus \{v_1,v_n\}$ and $J_n= V(P_{j_n})\setminus \{v_1,v_n\}$. Otherwise let $u$ be a vertex in $V(P_{j_1})$ such that removing $u$ splits $P_{j_1}$ into 2 subpaths each having at least $|S(v_1) \cap V(P_{j_1})|/2-1$ neighbors of $v$ in $S(v_1)$. Out of these 2 subpaths we let $P_u$ be the one with the most vertices in $S(v_n)$. If $|V(P_u)| \geq n/6-2$ then let $J_n=V(P_{u})\setminus \{v_1,v_n\}$ and $J_1=[n]\setminus (J_n \cup \{v_1,v_n\})$. Otherwise let $J_1=V(P_{j_1}) \setminus (V(P_u) \cup \{v_1,v_n\})$ and $J_n=[n]\setminus (J_1\cup \{v_1,v_n\})$. In both cases $J_i$  has size at least $n/6-4$, induces at most 2 subpaths  of $P$ and spans a set $S'(v_i)\subset S(v_i)$ of size 1000 such that the edges in $E(v)$ incident to $S(v_i)$ belong to the subpaths of $P$ induced by $J_i$ for $i\in \{1,2\}$. Let $B(j_i)$ be the set of vertices in $J_i$ that have a neighbor on $P$ not in $J_i$, $i=1,n$.

Let $\cP$ be the set of 1000 paths obtained by a P\'osa rotation with the inserting edge being $v_nv, v\in S'(v_n)$ and  $End$ the corresponding set of 1000 endpoints. Let  $S$ be a queue originally containing the elements in $End$. Throughout the algorithm for each element in both $S, End$ there will be a unique path in $\cP$.

While $|End| \leq n^{0.55}$ and $S\neq \emptyset$, let $v$ be the first vertex in $S$. Let $Rot$ be the elements in $u\in \{L_v(j):2q+1\leq j\leq 3q\}\cap J_n$ tha do not belong to $T_R$ for which the P\'osa rotation with the inserting edge being $vu$ results to a Hamilton path with the new endpoint  in $J_n \setminus (Q \cup End\cup B(J_n))$. Add $u$ to $Q$. If $|Rot|\geq 1000$ then add $Rot$ to $End$ and the elements of $Rot$ to the end of the queue $S$. In addition add the corresponding paths to $\cP$. 

Each Hamilton path in $P'\in \cP$  starts with $v_1$ and induces the same, at most 2, paths on $V(J_1)$. Thus it also induces 1 out of at most  4 possible orderings on $V(J_1)$ that we get by  traversing $P'$ starting from $v_1$. Let $\cP'$ be a subset of $\cP$ of size $|\cP|/4$ that agree on the ordering on $V(J_1)$ and let $End'$ be the set of the corresponding endpoints.

Now let $P'\in \cP'$ and $v_1,v'$ be its endpoints. Repeat the above procedure with $P'$ in place of $P$ to obtain a set of Hamilton paths $\cP''$ with a common endpoint $v'$ and let  $End''$ be the corresponding set of endpoints. Then 
examine the entries $\{l_v(i):v\in End'',i\in [2q+1,3q]$ until you find an entry $u$ that belongs to $End'$ or all of the corresponding entries have been examine. Add $End,End''$ to $Q$. 
If an edge $vu$ such that $v\in End'$, $u\in End''$ has been identified
then let $P_v\in \cP'$ and $P_u\in \cP''$ be the corresponding Hamilton paths.  Let $E(H_i)= E( P_v ) \triangle (E(P_u)\triangle E(P') ) \cup \{u,v\}$ and output $H_i=([n],E(H_i))$. If no such edge has been identified return FAILURE.

\begin{remark}
$E(P_v) \triangle (E(P_u)\triangle E(P'))$ is the edge set of a Hamilton path whose edges set is the union of the following edge sets: 
\begin{itemize}
    \item the edges of $P$ spanned by $B(J_1)\cup B(J_2)\cup\{v_1,v_n\}$.
    \item the edges of $P_v$ spanned by $(J_1 \cup \{v_1\}) \setminus B(J_1)$.
    \item the edges of $P_u$ spanned by $(J_2 \cup \{v_2\}) \setminus B(J_2)$.
\end{itemize}
\end{remark}
We give the analysis of the above algorithm in the following 2 Lemmas.

\begin{lemma}\label{lem:lists:determ}
During the $i^{th}$ iteration of Phase 3 at most $3qn^{0.55}$ queries are made and the size of $Q$ is increased by at most  $3n^{0.55}\leq n^{0.56}$ for $i\leq \ell$. In addition if the algorithm returns $H_i$ then for $e\in R\cap E(H_i)$ and $v\in e$ the edges in $E(v)\cap E(H_{i-1})$ belong to $E(H_i)$ and no vertex in $S(v)$ enters $Q$ during the $i^{th}$ iteration. 
\end{lemma}
\begin{proof}
Let $Q_i$ be the set of vertices that enters $Q$ at iteration $i$ of Phase 3. $Q_i\subset End \cup End''$ and by construction $|End|,|End''|\leq n^{0.55}+q$. Hence $|Q_i|\leq 2n^{0.55}+2q$. Thereafter, at iteration $i$, we query only entries from the lists of vertices in $Q_i$, and for each vertex  in $Q_i$ we make at most $q$ queries,  yielding at most $3qn^{0.55}$ queries during the $i^{th}$ iteration.

Now let $v\in e$ for some $e\in R\cap E(H_i)$. Observe that during the first $i$ iterations we never perform a P\'osa rotation where the pivot vertex belongs to $S(v)\cup T(v)$. Therefore no edge in $E(v)$ is deleted and no vertex in $S(v)$ ever enters one of $End, End''$; hence no vertex in $S(v)$ belongs to $Q$ by the end of the $i^{th}$ iteration of Phase 3.
\end{proof}

\begin{lemma}\label{lem:list:stoch}
At the end of the $i^{th}$ iteration of Phase 3, with probability $1-o(n^{-8})$ the following hold:
\begin{itemize}
    \item[(i)]  $|\cP|, |\cP'| \in [n^{0.55}, n^{0.55}+q]$.
    \item[(ii)] $H_i$ is outputted.
\end{itemize}
\end{lemma}

\begin{proof}
Recall that every set $Rot$ of size less than 1000 is discarded. Therefore for each set of 1000 consecutive elements from $S$ that are examined there exists at most 2 paths $P_1,P_2\in \cP$ such that each of those 1000 vertices is the new endpoint of a Hamilton path obtain from one of $P_1,P_2$ via a P\'osa rotation. Thus Lemma \ref{lem:exposerandomness} implies that among those 1000 vertices at most 600 belong to $D$. For each of those vertices $q$ queries are made each resulting to a new endpoint with probability, by Lemma \ref{lem:exposerandomness}, $\frac{|J_n|- 2|Q\cup  End \cup B(J_n)|-|V_T|}{n}>0.1$. Thus the probability that among those $400q>99\log n$ queries less than $2\log n$ potential new endpoints are identified is bounded by,
$$\Pr(Bin(99\log n,0.1)<2\log n)=o(n^{-9}).$$
Out of these endpoints at most $999\cdot 400$ are discarded. Thus with probability $1-o(n^{-8})$ among every 1000 consecutive elements from $S$ that are examined at least $\log n$ elements are added to each of $\cP, End$ and $S$. Therefore at the end of the while-loop we have that $S\neq \emptyset$ and $\cP\geq n^{0.55}$. As for every element of $S$ at most $q$ elements are queried from its list we have that $\cP\leq n^{0.55}+q$. Similarly $\cP'\leq n^{0.55}+q$.

Finally $H_i$ is not outputted only in the event that the desired $uv$ edge with $v\in End'$ and $u\in End''$ is not identified. From the above at least 0.4 portion of the vertices in each of $End',$ $End''$ do not belong to $D$. Thus, Lemma \ref{lem:exposerandomness}, implies that in the event $|End|,|End''|\geq n^{0.55}$ $H_i$ is not outputted with probability at most
$$\bigg(1-\frac{(0.4+o(1))|End'|}{n} \bigg)^{(0.4+o(1))|End''|}
\leq \bigg(1-n^{-0.46}\bigg)^{n^{0.54}} \leq e^{-n^{0.08}}=o(n^{-8}). 
$$

\end{proof}

Thus w.s.h.p. in Phase 3 we make at most $3\ell q n^{0.55}<n$ many queries yielding a total of at most $0.8\cdot 10^5n + n+n<10^5n$ queries through the 3 phase algorithm. In addition throughout Phase 3 at most $\ell \cdot q \cdot |Q| = O(n^{0.57}\log^2n)$ P\'osa rotations are performed yielding an $O(n)$ total running time of our 3-Phase Algorithm.

\section{Hamilton cycles in the adjacency matrix model}\label{sec:adjmatrix}
Let $G\sim G(n,p)$ for some $p=p(n)\geq 0$. In this section we assume that we are given the adjacency matrix $A_G$ of $G$. We consider that each of the entries of $A_G$ is blank and can be queried in $O(1)$ time. When we query the entry $A_G(i,j)$, or equivalently whether the edge $\{i,j\}$ belongs to $G$, we also reveal the entry $A_G(j,i)$ and we consider them to be the same entry. W.s.h.p. the algorithm that we implement will stop before  revealing $\big(1+\frac{1}{(\log\log n)^2}\big)n$ entries of $A_G$ that are equal to 1.
Therefore the Chernoff bound gives that  probability that it makes at least $\big(1+\frac{1}{\log\log n }\big)\frac{n}{p}$ queries is bounded by 
\begin{align*}
\Pr\bigg(Bin\bigg(\bigg(1+\frac{1}{\log\log n}\bigg )\frac{n}{p},p\bigg) < \bigg(1+\frac{1}{(\log\log n)^2}\bigg)n\bigg)=o(n^{-7}).
\end{align*}

For the rest of this section, modulo the last subsection, we  assume that $p\in \big[p^*,\frac{101\log n}{n}\big)$. We deal with the leftover cases at the end of this section. We also let $n_1=\frac{n}{\log^{0.6} n}$, $$V_1=[n-10n_1-1],\hspace{5mm} V_2=[n]\setminus V_1 \text{ and }V_{12}=\bigg[n-10n_1-1-\frac{n}{\log\log n},n-10n_1-1\bigg].$$

We prove Theorem \ref{thm:adjmatrix} by implementing a 2-phase algorithm. Its 1st phase is a greedy procedure that aims to build a small collection of vertex disjoint paths that cover $[n]$. It is divided into  a number of steps-algorithms. At the first one we greedily try to cover $[n]$ by a 2-matching that saturates all but $O(n_1)$ vertices. At the same time we distinguish 2 sets of vertices $D,V'$. Every vertex not in $D\cup V'$ will be saturated by $M$. In addition, every vertex in $V'$ is covered by $M$ and the subgraph of $G$ spanned by $V'$  will be distributed as a random graph on $V'$ with minimum degree 2. At the next couple of steps we will take care of $D$. At the final step we implement the Karp-Sipser algorithm \cite{karpsipser1981} to match almost all of the vertices in $V'$.

We exit this phase with a 2-matching $M$ and the sets $Burned$, $\cS$ and $V'$ that satisfy the following.
Each vertex not in $\cS \cup V'$ is saturated by $M$. Each vertex $i$ for which more than $10^{-5}n$ entries have been queried from row $i$ belongs to $Burned$. In addition, for each vertex  $v\in \cS\cup V'$ we also have put aside  sets $\cG_1(v), \cG_2(v), \cG_3(v)$ such that $\cG_1(v)$ contains at least $\log\log n-10^4$ neighbors of $v$,  $\cG_2(v)$ is the $M$-neighborhood of $G_1(v)$ and $\cG_1(v)\cup \cG_3(v)$ contains the $M$-neighborhood of $G_2(v)$.
Finally the sets $\{\cG_i(v)\}_{v\in \cS\cup V',i\in\{1,2,3\}}, Burned$ are pairwise disjoint.

The 2nd phase of our algorithm is almost identical to the Phase 3 algorithm given is Section \ref{sec:adjlists}.

\subsection{Phase 1: Covering the vertex set greedily}
Throughout the 1st phase $V'$ will be a subset of $V_2$. We would like to be able to certify that every vertex in $V'$ has at least 2 neighbors in $V'$ and still be able to treat those in the future as random. For this purpose it will be convenient to think that the algorithm is executed by two people, Alice and Bob. Alice, who plays the role of an oracle, goes ahead and reads a portion of the matrix. Thereafter she  keeps the information that Bob needs to access using appropriate functions which she occasionally updates. Then Bob executes the rest of algorithm. At any point he may read the functions that are kept by Alice and act accordingly. To see the usefulness of Alice and Bob consider the following example. We are given a row of size $n$ where each entry is an independent $Bernoulli(p)$ random variable. Alice goes first and reads the whole row and writes $c_2=1$ indicating that the corresponding row has at least two 1's. Then Bob reads $c_2=1$. At this point Alice knows the whole row and for her the row has become deterministic. However for Bob (and for us) the distribution of the row has changed but it has not become deterministic yet. Now for him the distribution is that the entries of the row are independent $Bernoulli(p)$  random variable conditioned on their sum being at least 2. He does not know which entries are equal to 1 but he is sure that if he starts querying the entries of the row one by one he will eventually find at least two 1's. 

Alice now reads the entries $A_G(i,j)$ with $j\in V_2$. $|V_2|= 10n_1$ implies that w.s.h.p. Alice reads at most $\frac{n}{\log^{0.1} n}$ 1's at this step. She then defines the function $c:V'\mapsto \{0,1,2\}$ given by $c(v)=0$ (respectively $=1,2$) if $v$ has 0 (resp 1 and at least 2) neighbors in $V'$. She then defines the sets $C_i$, $i=1,2,3$ by $C_i=\{v\in V':c(v)=i\}$. Henceforward if $V'$ changes and an algorithm asks to update the function $c()$ or the sets $C_i,i=1,2,3$ then Alice proceeds and performs this step. Bob, who is aware of the updates, performs the rest of the steps of the algorithm. We will mark the steps performed by Alice using (Alice). 

A procedure that Alice will implement is to identify the 2-core (i.e. the maximal subgraph of minimum degree 2) of a subgraph of $G[V_2]$. For that she implements the following procedure for $V''\subset V_2$.

\begin{algorithm}[H]
\caption{Identify2Core($V''$)}
\begin{algorithmic}[1]
\While{ there exist $v\in V''$ with 0 or 1 neighbors in $V''$}
\\ \hspace{5mm} Remove $v$ from $V''$
\EndWhile
\\Return $V''$.
\end{algorithmic}
\end{algorithm}

At GreedyCover, the algorithm given below, $M$ is a 2-matching, $S$ is the set of $M$-saturated vertices and $C$ is the set of vertices in $[n]$ that are covered but are not saturated by $M$. Throughout its execution  $M$  induces a set of paths. If $v$ is the endpoint of such a path $P$ then $e(v)$ is the other endpoint of $P$  and the algorithm will never query $A_G(v,e(v))$. Thus $M$ does not induce a cycle at any point of the execution of GreedyCover. In addition no entry $A_G(i,j)$ with $i,j\in V_2$ is queried. The functions $count(\cdot), count'(\cdot)$ ensure that initially, not too many entries from any row are queried. At the end of CoverGreedy, the set $D$ will be the union of the sets $D_0$, $D_1$, $D_2$ and will consists of (potentially) problematic vertices. $D_0$ will consists of vertices in $V_{2}$ that do not belong to the 2-core of some subgraph of $G$ with vertices in $V_2$. $D_1$ will consists of non $M$-saturated vertices in $V_1$ and non $M$-covered vertices in $V_2$. $D_2$ will consist of vertices in $V_2$ with potentially small degree in $G$. After forming $M$ we put aside a sets of vertices $Extend$, containing to edges with an endpoint in $D\cup V_2$. We will use these edges  to reroute the paths induce by $M$ and to join the paths into a Hamilton cycle later on. Finally the set $Burned$ will contain all the vertices $v$ for which we have queried at least $10^{-5}n$ entries from the corresponding row.

For $v\in [n]$ and a 2-matching $M$ of $G$ we let $N_M(v)$ be the neighbors of $v$ in $G$ that are connected to $v$ via an edge in $M$ and $d_M(v)=|N_M(v)|$. Similarly for $U\subseteq [n]$ we let $N_M(U)$ be the set of vertices in $[n]\setminus U$ that are adjacent to $U$ via an edge in $M$. Finally, we refer to the graph-distance induced by the graph $([n],M)$ as $M$-distance.
 
\begin{algorithm}[H]
\caption{GreedyCover}
\begin{algorithmic}[1]
\\$D_0=D_1=D_2=D=S=C=M=Burned=Extend=\emptyset$. 
\\ For $i\in [n]$ set $e(i)=i$ and $count(i)=count'(i)=0$
\For{ $i \in V_2$}
\While{ $i\notin C$
 and $count'(i)<n_1$ }
\\\hspace{10mm} Let  $j= \min\{ k\in V_{12}\setminus (S \cup \{e(i)\}):count(k)<2n_1 \text{ and }A_G(i,j) \text{ has not been queried}\}$.
\\\hspace{10mm} Query $A_G(i,j)$ and update $count(j)=count(j)+1$, $count'(i)=count'(i)+1$.
\\\hspace{10mm}  If $A_G(i,j)=1$ then add $\{i,j\}$ to $M$ and set $S=S\cup (\{i,j\} \cap C)$, $C=\big(C\cup \{i,j\}\big) \setminus S$, $e(e(i))=e(j)$ and $e(e(j))=e(i)$.
\EndWhile
\EndFor
\For{ $i \in V_1$}
\While{ $i\notin S$
 and $count'(i)<n_1$ }
\\\hspace{10mm} Let  $j= \min\{ k\in [n]\setminus (S \cup \{e(i)\}):count(k)<2n_1 \text{ and }A_G(i,j) \text{ has not been queried}\}$.
\\\hspace{10mm} Query $A_G(i,j)$ and  update $count(j)=count(j)+1$, $count'(i)=count'(i)+1$.
\\\hspace{10mm}  If $A_G(i,j)=1$ then add $\{i,j\}$ to $M$ and set $S=S\cup (\{i,j\} \cap C)$, $C=\big(C\cup \{i,j\}\big) \setminus S$, $e(e(i))=e(j)$ and $e(e(j))=e(i)$.
\EndWhile
\EndFor
\For{ $v\in  V_2$}
\\\hspace{5mm} Query the entries from $v$ to $V_{12}\setminus Extend$ one by one until $\log\log n$ edges are detected or all of the corresponding entries have been queried and let $\mathcal{G}(v)$ be the corresponding set of neighbors of $v$. 
\If{ $|\mathcal{G}(v)|<\log\log n$}
\\\hspace{10mm} Add $v$ to $D_2$.
\Else
\\\hspace{10mm} Add to $\cG(v)$ all the vertices within $M$-distance 4  from $\mathcal{G}(v)$.
\\\hspace{10mm} Add $\cG(v)$ to $Extend$ 
\EndIf
\EndFor
\\ Let $D_1=(V_1\setminus S) \cup (V_2 \setminus C)$
\\(Alice) $V'=$Identify2Core($V_2\setminus (D_2\cup D_1))$.
\\ Set $D_0= (V_2\setminus (D_2\cup D_1)) \setminus V'$
and $D=D_0\cup D_1\cup D_2$.
\For{ $v\in  D$}
\\\hspace{5mm} Query the entries from $v$ to $V_1 \setminus (V_{12}\cup D\cup Extend)$ until $0.01\log n$ edges are detected or all of the corresponding entries have been queried and let $E(v)$ be the corresponding set of neighbors of $v$. 
\\\hspace{5mm} If $|E(v)|=0.01\log n$ then add to $E(v)$ all the vertices within $M$-distance 4 from $E(v)$.
\\\hspace{5mm} Add $E(v)$ to $Extend$.
\EndFor
\\If $|D| > \frac{n}{\log^9n}$ then return FAILURE.  
\\ Set $Burned=D$.
\\Return $V',D, S, C,Burned, Extend, \{\cG(v)\}_{v\in V_2}, \{E(v)\}_{v\in D}, M$ 
\end{algorithmic}
\end{algorithm} 
\begin{observation}\label{obs:Gnm2}
Upon termination of GreedyCover, given the number of edges spanned by $G[V']$, say $m$, $G[V']$ is distributed as a random graph on $V'$ with $m$ edges and minimum degree 2.
\end{observation}

\begin{lemma}\label{lem:greedycover}
W.s.h.p. $|D| \leq \frac{n}{\log^9 n}$ and $|V'|\geq 5n_1$. In addition, w.s.h.p. GreedyCover does not return FAILURE and it terminates before it reads $\big(1+\frac{1}{\log^{0.3} n}\big)n$ 1's.
\end{lemma}

\begin{proof}
A vertex $i\in V_1$ (respectively in $V_2$) does not belong to $S$ (resp. $C$) if at the end of GreedyCover $count'(i) < n_1$. Therefore,
\begin{align*}
    \Pr&\bigg(|D_1| \geq \frac{n}{3\log^9n}\bigg)
    \leq \binom{n}{\frac{n}{3\log^9 n}}
    \big[(1-p)^{n_1-1}+n_1p (1-p)^{n_1-2}\big]^{\frac{n}{3\log^9 n}}
\\& \leq \bigg( (e\log^9 n)\cdot 6n_1p \cdot e^{-pn_1}    \bigg)^{\frac{n}{3\log^9 n}}
=o(n^{-7}).
\end{align*}
The last equality holds as for $p\geq p^*$ we have that $n_1p\geq0.99 \log^{0.4}n$.
On the other hand a vertex $v\in V_2$ belong to $D_2$ only if at line 18 fewer than $\log\log n$ neighbors are detected. In the event $v\in D_2$ at line 18 at least $\frac{n}{\log\log n}- O(|V_2|\log\log n)\geq \frac{n}{2\log\log n}$ unqueried entries are queried. As $|V_2|\leq n$, 
\begin{align*}
    \Pr&\bigg(|D_2| \geq \frac{n}{3\log^9n}\bigg)
    \leq \binom{n}{\frac{n}{3\log^9 n}}
\bigg[\Pr\bigg(Bin\bigg(\frac{n}{3\log\log n},p\bigg)\leq \log\log n\bigg)\bigg]^{\frac{n}{3\log^9 n}}
\\&\leq \bigg(3e\log^9n \Pr\bigg(Bin\bigg(\frac{n}{2\log\log n},p\bigg)\leq \log\log n\bigg) \bigg)^{\frac{n}{2\log^9 n}}
=o(n^{-7}).
\end{align*}
At the last equality we used the Chernoff bound. In the event that $|D_1\cup D_2|\leq \frac{2n}{3\log^9 n}$ and $|D_0|> \frac{n}{3\log^9 n}$ by considering the execution of Identify2Core we have that there exists disjoint sets $U,W, Z\subseteq V_2$ such that $U=D_1\cup D_2$ , $|U|\leq \frac{2n}{3\log^9 n}$, $Z=N_{G[V_2]}(W)$, $W\subseteq D_0$, $|W|= \frac{n}{3\log^9 n}$ and every vertex in $W$ has at most 1 neighbor in $V_2$, thus $|Z|\leq |W|$. Hence,
\begin{align*}
    \Pr\bigg(|D_0|>\frac{n}{3\log^9 n} \bigg| |D_1\cup D_2| \leq \frac{2n}{3\log^9 n}\bigg)\leq  \binom{10n_1}{ \frac{2n}{3\log^9 n}}\binom{10n_1}{ \frac{n}{3\log^9 n}}\binom{10n_1}{ \frac{n}{3\log^9 n}} (1-p)^{(1+o(1)) \frac{|V_2|n}{3\log^9 n}}
    \\ \leq \bigg(\log^{40}n\cdot e^{-(1+o(1))p|V_2|} \bigg)^{\frac{n}{3\log^9 n}} \leq \bigg(\log^{40}n\cdot e^{-\log^{0.4} n} \bigg)^{\frac{n}{3\log^9 n}}=o(n^{-7}).
\end{align*}

Now, if $i\in V_1$ is the first vertex to be matched with a vertex in $V_2$ (at the for-loop at line 10) then as $count(i)<2n_1$ it must  be the case that when $i$ is matched there exists at most $2n_1-1$ unsaturated vertices between $i$ and $n-10n_1$ that are queried after $i$ and may be matched to vertices in $V_2$ afterwards. Thus, $|S \cap V_2|\leq 4n_1$ and w.s.h.p. $|V'|\geq |V_2\setminus S|-|D| \geq 5n_1$.

Finally, GreedyCover identifies at most
$\Delta(G)$ edges from $s$ to $[n]$ for every $s\in D$, at most 2 edges from $v$ to $V_1\setminus D$ for every $v\in V_1\setminus D$, at most $\log\log n$ edges from $v$ to $V_1\setminus D$ for every $v\in V_2\setminus D$ and no edge spanned by $V_2$. 
At Lemma \ref{lem:satdanger:basicprop} we show that w.s.h.p. $\Delta(G)\leq 200\log n$. Thus w.s.h.p. GreedyCover reads at most $\frac{200n}{\log^8 n}+n+\frac{n\log\log n}{\log^{0.6}n}\leq \big(1+\frac{1}{\log^{0.3}n}\big)n$ 1's.

\end{proof}

\subsection{Phase 1: Saturating the troublemakers}

In this step we slightly alter $M$ such that it saturates a set $S'$, it covers $V'$ and $|[n]\setminus (S' \cup V')| =O(n^{0.1})$. We do so by implementing first SatDangerous and then SatTrouble. Let $Dangerous$ be the set of vertices $v\in D$ with $|E(v)|<0.01\log n$. SatDangerous aims to alter $M$ such that all the vertices in $Dangerous$  are saturated. It explores their neighborhood and finds a set of vertex disjoint paths that cover them which then uses to alter $M$. Later we show that w.s.h.p. $|Dangerous|\leq n^{0.1}$, thus the overall change  of $M$ at this step is minor. The set of paths are carefully chosen such that any new unsaturated vertices have significantly large degree. Afterwards SatTrouble takes care of the rest of the vertices in $D$.

The sets $V',S,C,D,\{E(v)\}_{v\in D}, \{\cG(v)\}_{v\in V'},Extend,Burned,$ and  the 2-matching $M$  $A_G$ are given by GreedyCover.

\begin{algorithm}[H]
\caption{SatDangerous($M,A_G,D, V',Burned, Extend, Dangerous \{E(v)\}_{v\in D}, \{\cG(v)\}_{v\in V'}$)}
\begin{algorithmic}[1]
\\ Define the function $q:V\mapsto \{0,1\}$ by $q(v)=0$ for $v\in V$.
\\ Define the graph $F=(V(F),E(F))$ by  $V(F)=Dangerous$ and $E(F)=\emptyset$.  
\While{ there exists $v\in V(F)$ with $q(v)=0$ such that either $v$ or an $M$-neighbor of $v$ is adjacent to $Dangerous$ } 
\\\hspace{5mm} Let $j=1$ and $S(v)=\emptyset$.
\While{$j\leq n$ and $|S(v)\setminus (V' \cup D  \cup Extend\cup Burned)|<0.01\log n$ }
\\\hspace{10mm} Reveal $A_G(v,j)$. If $A_G(v,j)=1$ then add $j$ to $S(v)$. 
Set $j=j+1$.
\EndWhile
\\\hspace{5mm} If $|S(v)\setminus (V' \cup D  \cup Extend\cup Burned)|<0.01\log n$ then add $v$ to $Dangerous$.
\\\hspace{5mm} Set $q(v)=1$. Add to $S(v)$ all of its $M$-neighbors. 
\\\hspace{5mm} Add $S(v)$ to $V(F)$ and the corresponding edges to $E(F)$.
\EndWhile
\\ If $|Dangerous'|>n^{0.1}$ or $F$ has a component with more than 9 vertices in $Dangerous$ then return FAILURE.
\\ For each component of $F$, by examining all the possible combinations of edges incident to Dangerous, search for a set of vertex disjoint paths $\cP$ that satisfies the following:
\begin{itemize}
    \item All the vertices in $Dangerous$ lie in the interior of some path in $\cP$.
    \item The paths in $\cP$ cover in total at most $3|Dangerous|$ many vertices.
    \item The paths in $\cP$ are spanned by  $Dangerous\cup N(Dangerous)$.
\end{itemize}
\If{ not such a set $\cP$ exists }
\\ \hspace{5mm} Return ``$G$ is not Hamiltonian", $Dangerous$.
\EndIf
\\  Let $V_P$ be the set of vertices covered by some path in $\cP$.
\\  Remove from $M$ all the edges incident to $V_P$.
\\  Let $Special'$ be the set of vertices whose $M$-degree has just been decreased.
\\  For $v\in Special'$ let $F(v)$ be the set of neighbors of $v$ in $V(F)$ not in $V'\cup D\cup Extend\cup Burned$.
\\  Add to $Burned$ every vertex $v\in V(F)$ that satisfies $q(v)=1$.
\For{ $v\in Special'$ }
\\ \hspace{5mm} Let $\cG(v)=F(v)$.
\\ \hspace{5mm} Add to $\cG(v)$ and all vertices within $M$-distance 4 from $\cG(v)$. 
\\ \hspace{5mm} Set $\cG'(v)=\cG(v)\setminus (Extend \cup Burned)$.
\\ \hspace{5mm} Add to $\cG'(v)$ to $Extend$.
\EndFor
\\(Alice) $V''=Identify2Core(V'\setminus Special')$.
\\ Set $Special=Special' \cup (V'\setminus V'')$ and $V'=V''\setminus (Dangerous \cup Special)$.
\\ Set $D=D\setminus (Dangerous \cup Special)$.
\\  Set  $E'(v)=E(v)\setminus (Extend \cup Burned \cup Special)$ for $v\in D$.
\\  Set
$\cG'(v)= \cG'(v) \setminus (Burned \cup Extend \cup Special)$ for $v\in V'$.
\\  If there exist $v\in Special'$ such that $|F(v)|<0.01\log n$ or $|\cG(v)\setminus \cG'(v)|\geq 10^4$ then return FAILURE.
\\  If there exist $v\in D$ such that $|E(v)\setminus E'(v)| > 10^4$ then return FAILURE.
\\  If there exist $v\in V'\cup (Special\setminus Special')$ such that $|\cG(v)\setminus \cG'(v)| > 10^4$ then return FAILURE.
\\ Return $M$, $D$, $V'$, $Special$, $Extend$, $\{E'(v)\}_{v\in D}$, $\{\cG'(v)\}_{v\in V'\cup Special}$ and $Burned$. 
\end{algorithmic}
\end{algorithm} 

\begin{observation}\label{rem:failure1:repeat}
If SatDangerous returns ``$G$ is not Hamiltonian" then this statement is true and $Dangerous$ is a certificate of this fact. Indeed as $N_F(Dangerous)=N_G(Dangerous)$ if  $H$ is a Hamilton cycle of $G$ then we have that the edges in $E(H)$ incident to $Dangerous$ induce a set of paths $\cP$ with the desired properties. 
\end{observation}
For $v\in D$  let $E_1'(v)$ be the neighbors of $v$ in $E'(v)$ that have 6 vertices within $M$-distance 3 in $E'(v)$. In addition let $E_2'(v)=N_M(E_1'(v))$ and $E_3'(v)=N_M(E_1'(v))\setminus E_1'(v)$. Similarly define the sets $\cG_1'(v),\cG_2'(v),\cG_3'(v)$ for $v\in V'\cup Special$.
\begin{observation}\label{rem:success:sattrouble}
If SatDangerous does not return FAILURE then,  $|E_1'(v)|\geq 0.01\log n-10^4$ and $|\cG_1'(u)|\geq \log\log n-10^4$ for $v\in D$ and $u\in V'\cup Special$ respectively. In addition the sets $\{E_i'(v)\}_{v\in D,i=1,2,3}$, $\{\cG_i'(v)\}_{v\in V'\cup Special, i=1,2,3}$ are pairwise disjoint, do not intersect $Burned$ and are subsets of $Extend$.
\end{observation}

To analyze SatDangerous we will make use of the following lemma. Its proof is given at Appendix \ref{app:matrix}. 

\begin{lemma}\label{lem:satdanger:basicprop}
W.s.h.p. the following hold:
\begin{itemize}
\item [(i)] $\Delta(G)\leq 200\log n$.
\item[(ii)] There does not exist a connected subgraph of $G$ spans $k$ vertices and more than $k+8$ edges for $k\leq 10^5$.
\item[(iii)] $|Dangerous|\leq n^{0.1}$.
\item[(iv)] No connected subgraph of $G$ spans $k\leq 1000$ vertices out of which more than 9 vertices belong to $Dangerous$.
\end{itemize}
\end{lemma}

\begin{lemma}\label{lem:satdangerous}
W.s.h.p. SatDangerous does not return FAILURE and $|Special'|\leq 6n^{0.1}$.
\end{lemma}
\begin{proof}
If SatDangerous returns FAILURE then either (I) $|Dangerous|>n^{0.1}$ or (II) $F$ has a component with at least 9 vertices from $Dangerous$ or (III) FAILURE is return at one of the lines 33-35. 

Parts (iii), (iv) of Lemma  \ref{lem:satdanger:basicprop} imply that w.s.h.p. (I),(II) do not occur. Let $\cE_{33}$ be the event that SatDangerous returns FAILURE at line 33 for some vertex $v$ in $Special'$. As every vertex $v\in Special'$ satisfies the condition at line 8 we have that $|F(v)|\geq 0.01\log n$  for $v \in Special'$. 

Now observe that every vertex in $Special'$ is at distance at most 2 form $Dangerous$. Thus, for $v\in Special'$ and for each $u\in \cG(v)\setminus \cG'(v)$  there exists a path of length at most 7 from $v$ to $Special'\cup V_P$, hence a path of length at most 9 from $v$ to $Dangerous$ that covers the vertex $u$. Therefore if $\cE_{33}$ occurs then either there exists $v\in V\setminus Dangerous$ and more than 9 vertices in $Dangerous$ within distance at most 10 from $v$ in $G$ or there exists sets $S,T,U$ of size 1, at most 9  and $1000<10^4$ in $V(G)$ such that $U$ is covered by $S$-$T$ paths of length at most 9. In the second case $G$ spans a subgraph on $k$ vertices with at least $k+100$ edges for some $k\leq 10^4$. Lemma \ref{lem:satdanger:basicprop} implies that w.s.h.p. none of the 2 cases occur and therefore $\Pr(\cE_{33})=o(n^{-7})$. 

The proof of the statement that w.s.h.p. SatDangerous does not return FAILURE at lines 34 and 35 follows in a similar fashion.

Now to bound the size of $Special'$ observe that for every vertex in $V_P$ we remove at most 2 edges from $M$ at line 18 and therefore add at most 2 vertices to $Special'$. Thus w.s.h.p. 
$$|Special'|\leq 2|V_P|\leq 6|Dangerous|\leq 6n^{0.1}.$$
\end{proof}
We proof the following lemma in Appendix \ref{sec:app:KS}.
\begin{lemma}\label{lem:special:gnm2}
W.s.h.p. $|Special|\leq n^{0.3}$.
\end{lemma}
Lemmas \ref{lem:greedycover},\ref{lem:satdangerous} and \ref{lem:special:gnm2} imply that after the execution of SatDangerous we have that w.s.h.p. $|V'|\geq 4n_1$.

For   $v\in V$ we let $end(v)$ be the set of endpoints of the path that is induced by $M$ and covers $v$. We also let $Forbidden=\cup_{v\in D}end(v)$. Lines 18-28 of SatTrouble also appear in the description of Karp-Sipser. 

\begin{algorithm}[H]
\caption{SatTrouble($M,A_G,D,V',\{E'(v)\}_{v\in D}, \{end(v)\}_{v\in [n]}, Forbidden$)}
\begin{algorithmic}[1]
\\ Set $Special_0=Special_1=\emptyset$.
\While{  $D\neq \emptyset$  }
\\\hspace{5mm}  Let $v \in D$. 
\While{$E'(v)\neq \emptyset$ and $d_M(v)\neq 0$.} 
\\\hspace{10mm} Let $u \in E'(v)$. Query the edges from $u\cup N_M(u)$ to $V'\setminus Forbidden$. 
\\\hspace{10mm} Add $u\cup N_M(u)$ to $Burned$. 
\If{ each vertex in $u\cup N_M(u)$ is adjacent to at least 7 vertices in $V'\setminus Forbidden$}
\\\hspace{15mm} From $M$ remove the edges from $u$ to $N_M(u)$.
\\\hspace{15mm} Add $vu$ to $M$ and $end(u)$ to Forbidden. 
\For {$w\in\{u\} \cup N_M(u)$}
\\ \hspace{20mm} Match $w$ to $w'\in V'\setminus Forbidden$, add $ww'$ to $M$, add $end(w')$ to Forbidden.
\\ \hspace{20mm}  Remove $w'$ from $V'$.
\EndFor
\EndIf
\EndWhile
\\\hspace{5mm} If $d_M(v)\neq 0$ then return FAILURE.
\\\hspace{5mm} Update the sets $C_0,C_1,C_2$.
\While{ $C_1 \neq \emptyset$ }
\\ \hspace{10mm} Let $v\in C_1$. Let $vu$ be the unique edge incident to $v$ in $G[V']$.
\If{ $u= end(v)$ }
\\ \hspace{15mm} Add $u,v$ to $Special_1$.
\Else
\\ \hspace{15mm} Add $uv$ to $M$ and remove $v,u$ from $V'$.
\\ \hspace{15mm} Update, with $v_1=end(v)$ and $u_1=end(u)$, $end(v_1)=u_1$ and $end(u_1)=v_1$.
\EndIf
\\ \hspace{10mm} Update the function $c()$ and sets $C_0,C_1,C_2$.
\\ \hspace{10mm} Remove $C_0$ from $V'$ and add it to $Special_0$.
\EndWhile
\\ Return $M$, $Special_0,Special_1,V',Burned$.
\EndWhile
\end{algorithmic}
\end{algorithm}

\begin{lemma}\label{lem:Q}
W.s.h.p. no more than $\frac{n}{\log^8 n}$ vertices are removed from $V'$ during the execution of SatTrouble. 
\end{lemma}
\begin{proof}
Let $U$ and $W$ be the set of vertices in $V'$ that are removed during the execution of SatTrouble because of lines 4-15 and 18-28 respectively. Then $|U|\leq 6|D|$, hence by Lemma \ref{lem:greedycover} we have that $|U|\leq \frac{6n}{\log^8 n}$. Thereafter for every vertex $w \in W$ at the moment when $w$ is removed from $V'$ either $w$ has at most 1 neighbor in $V'$ or it is the unique neighbor in $V'$ of some vertex $w'\in V'$.  By considering the moment that the $\bfrac{n}{\log^8 n}$th vertex is removed from  $V'$ we have that $V'$ can be partitioned into 4 set $U,W_1,W_2,Z$ such that 
$|U|\leq \frac{6n}{\log^9n}$,  $|W_1|\geq \frac{0.5n}{\log^8 n}-\frac{2.5n}{\log^9n}$ and $|W_2|\leq |W_1|$ such that no vertex in $W_1$ is adjacent to $Z$.

Therefore the probability that more than $\frac{n}{\log^8 n}$ are removed from $V'$  during the execution of SatTrouble
is bounded by,
\begin{align*}
    &\binom{|V'|}{\frac{6n}{\log^9 n}}
    \binom{|V'|}{\frac{0.5n}{\log^8 n}}     \binom{|V'|}{\frac{0.5n}{\log^8 n}}
    (1-p)^{|V'| \cdot \frac{0.49n}{\log^8 n}}    \\&\leq \bigg((\log n)^{20}\cdot e^{-0.49p|V'|} \bigg)^\frac{n}{\log^8 n}
    \leq \bigg((\log n)^{20}\cdot e^{-0.4\log^{0.4}n} \bigg)^\frac{n}{\log^8 n}
 =o(n^{-7}).
\end{align*}
At the last inequality we used that w.s.h.p. at the beginning of SatTrouble $|V'|\geq 4n_1$.
\end{proof}

\begin{lemma}\label{lem:satterm}
W.s.h.p. SatTrouble does not return FAILURE.
\end{lemma}
\begin{proof}
SatTrouble returns FAILURE while a vertex $v\in D$ is examine at the while-loop at line 4 if initially there do not exists 2 vertices $w,w'\in E_1'(v)$ such that every vertex in $W=\{w,w'\}\cup N_M(\{w,w'\})$ has at least 14 neighbors in $V'$. Observation \ref{obs:3} implies that $|E_1'(v)|\geq 0.01\log n-10^4$ and for $w\in E_1'(v)$ and $u\in \{w\}\cup N_M(w)$ no edge from $u$ to $V'$ has been queried so far. In addition Lemma \ref{lem:Q} implies that w.s.h.p. $|V'|\geq 4n_1-\frac{n}{\log^8n}\geq 3n_1$. Thus, the probability that a vertex $w\in E_1'(v)\cup E_2'(v)$ has fewer than 14 neighbors in $V'$ is less than 
$$\Pr(Bin(3n_1,p)\leq 14) \leq e^{-1.5n_1p} \leq e^{-\log^{0.4} n}$$
and the probability that SatTrouble returns FAILURE  is at most
\begin{align*}
    n \bigg(\Pr(Bin(0.01\log n-10^4,3\cdot e^{-log^{0.4} n})\geq 0.01\log n-10^4-2 \bigg)=o(n^{-7}).
\end{align*}
\end{proof}
\subsection{Phase1: The Karp-Sipser algorithm}
We reach step 4 with a 2-matching $M$ and sets $V',Special,Special_0,Special_1\subset V$ such that $M$ saturates $[n]\setminus (V' \cup Special\cup Special_0\cup Special_1)$ and covers $V'$. In addition $M$ induces a set of vertex disjoint paths on $[n]$. Furthermore, no edge spanned by $V'$ has been queried (by Bob). We now implement the Karp-Sipser algorithm as stated below to extend $M$ to a 2-matching of size at least $n-2n^{0.41}$ that spans no cycle. 

Recall, Alice keeps and updates the sets $C_0,C_1$ and $C_2$. Initially $C_0=C_1=\emptyset$ and $C_2=V'$.

\begin{algorithm}[H]
\caption{Karp-Sipser($M,V',\{end_v\}_{v\in V'}$, $Special_0,Special_1$)}
\begin{algorithmic}[1]
\While{ $|V'|\leq n^{0.41}$ }:
\If{ $C_1\neq \emptyset$ }
\\ \hspace{10mm} Let $v\in C_1$. Let $vu$ be the unique edge incident to $v$ in $G[V']$.
\Else
\\ \hspace{10mm}
Let $v\in C_{2}$. Reveal the neighbors of $v$ in $V'$, and let $u$ be one of them other than $end(v)$.
\EndIf
\If{ $u= end(v)$ }
\\ \hspace{15mm} Add both $v,u$ to $Special_1$.
\Else
\\ \hspace{15mm} Add $uv$ to $M$.
\EndIf
\\ \hspace{5mm} Remove $v,u$ from $V'$.
\\ \hspace{5mm} (Alice) Update $c()$ and the sets $C_0,C_1,C_2$. \\ \hspace{5mm} Add $C_0$ to $Special_0$ and remove it from $V'$.
\\ \hspace{5mm} Update, with $v_1=end(v)$ and $u_1=end(u)$, $end(v_1)=u_1$ and $end(u_1)=v_1$.
\EndWhile
\\ Output $M$, $Special_0, Special_1$.
\end{algorithmic}
\end{algorithm}

\begin{lemma}\label{lem:karp2}
W.s.h.p. $|Special_1|\leq n^{0.41}/3$.
\end{lemma}

\begin{proof}
Let $uv$ be some edge detected either by SatTrouble at lines 18-28 or by Karp-Sipser. At the time when $v\in C_1$ is chosen Bob has not revealed any of the edges that are spanned by $V'$ (as it was defined at that point during the execution of the corresponding algorithm). Hence, given the number of edges spanned by $V'$, say $m$, the current graph is uniform over all the graphs on $V'$ with $m$ edge in which every vertex in $C_1$ has degree 1 and every other vertex has degree at least 2.  $v,u$ are added to $Special_1$ only if $u$ is the unique neighbor of $v$ in $V'$. 
The probability of this event occurring is at most $1/|V'|$ if $u\in C_1$ and at most $\frac{1}{|V'|-|C_1|}\leq \frac{2}{|V'|}$ (see Section \ref{sec:app:KS} of the appendix for further explanation). Therefore the number of vertices in $Special_1$ is dominated by twice the sum of $n$ variables $X_2,X_3,...,X_n$ where $X_i$ is a $Bernoulli(2/i)$ random variable.
Thus the probability that $|Special_1|\geq n^{0.41}/3$ is bounded by the probability there exists $i\leq \log_2 n-1 $ such that $\sum_{j=2^i+1}^{2^{i+1}}X_i > \log^2 n$.
Hence,
\begin{align*}
\Pr(|Special_1|\geq n^{0.41}/3)&\leq \sum_{i=0}^{\log_2n-1} \Pr(Bin(2^i,2^{-i+1})\geq \log^2 n)   
\\& \leq \sum_{i=0}^{\log_2n-1} \binom{2^i}{\log^2 n} \bfrac{2}{2^i}^{\log^2 n} 
\\& \leq \sum_{i=0}^{\log_2n-1} \bfrac{2e}{0.5\log n} ^{\log^2 n}=o(n^{-7}).
\end{align*}
\end{proof}

By slightly modifying the arguments given in 
\cite{anastos2021k}, \cite{anastos2020} and   \cite{Hamd3}, one can prove Lemma \ref{lem:karp}. For sake of completeness we display the key elements of the corresponding arguments at Appendix \ref{sec:app:KS}.

\begin{lemma}\label{lem:karp}
W.s.h.p. $|Special_0|\leq n^{0.41}/3$.
\end{lemma}

\subsection{Phase 2: Transforming the path packing into a Hamilton cycle}
Let $\cS=Special\cup Special_1\cup Special_2 \cup V'$. Lemmas \ref{lem:satdangerous}, \ref{lem:karp2} and \ref{lem:karp} imply that $|\cS|\leq 2n^{0.41}$. We reach Phase 2 with a 2-matching $M$ such that $n-|M|\leq |\cS| \leq 2n^{0.41}$ and sets $\cS,\{\cG_i'(v)\}_{v\in \cS,i=1,2,3}$, $Burned$ and $Extend$. Each vertex $v$ such that more than $10^{-5}n$ edges incident to $v$ have been queried belongs to $Burned$.
Thereafter $\cG_1'(v)$  contains at least $\log\log n-10^4$ neighbors of $v$, $\cG_2'(v)=N_M(\cG_1'(v))$ and $\cG_3'(v)=N_M(\cG_2'(v))\setminus \cG_1'(v)$. In addition the sets $\{\cG_i'(v)\}_{v\in \cS,i=1,2,3},\cS$ and $Burned$ are pairwise disjoint,  $\bigcup_{v\in \cS,i=1,2,3}
N_M(\cG_i'(v)) \subseteq Extend$ and $|Extend|=o(n)$.
Finally $M$ does not span a cycle and saturates every vertex in $[n]\setminus \cS$. Thus w.s.h.p. it induces at most $2n^{0.41}$ path components whose endpoints belong to $\cS$. Here we allow for paths of length 0 which correspond to vertices  not covered by $M$.

Let $P_1,P_2,...,P_k$ be those paths. Henceforward we implement an algorithm almost identical to the one given for Phase 3 in Section 2. There are only 2 major differences. The first one is that we now allow for paths of length 0. To cope with this, after introducing a set of edges $R$ and joining the paths to a Hamilton cycle $H$, for each path-vertex $v$ we split the set $\cG_1'(v)$ into two sets of the same size; let them be $\cG_{11}'(v)$ and $\cG_{12}'(v)$. Let $e_1,e_2$ be the edges in $R$ incident to $v$ in $H$. We use $\cG_{11}'(v)$ (and $\cG_{12}'(v)$ respectively) at the iteration that starts with removing from the current Hamilton cycle the edge $e_1$ ($e_2$ resp.).

The second difference is that at Section 2, given some endpoint $v$  we revealed some entries of $L_v$, we kept those that lie in $J_n$ and resulted to a new endpoint $u$ for which we did not query entries from $L_u$ at this phase and then did the corresponding P\'osa rotations. We then repeated this step until $End$ becomes of size $n^{0.55}$. Now, given $v\in End$ we query the edges from $v$ to the vertices in $J_n$ that will result to a new endpoint not in $Extend\cup Burned$. We then add $v$ to $Burned$. In total we  perform at most $2\cdot 2n^{0.41}$ sequences of P\'osa rotations, each time producing $n^{0.55}$ many endpoints which are then added to $Burned$. Therefore throughout the algorithm $|Burned|=o(n)$.
For each endpoint we query  $|J_n|-o(n)\geq 0.1n$ edges each resulting to a new endpoint independently with probability $p$. At this point it is important to point out that the sets $\{\cG_i'(v)\}_{v\in \cS,i=1,2,3}$ and $Burned$ are pairwise disjoint and are also subsets of $Extend$. Therefore for $v\in \cS$ and $u\in \cG_2'(v)$ the vertex $u$ may appear as an endpoint of a Hamilton path in $G\cup R$ only in one case. That is the case   where we start with a Hamilton path $P$ which has $v$ as an endpoint and the pivot vertex is the common neighbor of $v,u$ that lie in $\cG_1'(v)$.

Finally w.s.h.p. the number of 1's that our algorithm read in total is at most
$$\bigg(1+\frac{1}{\log^{0.3}n} \bigg)n+(|D|+|Dangerous|)\log^3n+ 2p|V_2|^2  +n^{0.97}\leq \bigg(1+\frac{1}{(\log\log n)^2}\bigg)n.$$ 
\subsection{Small or large p}
If $p\geq n^{-1/3}$ then we start by  running lines 1-25 of Greedy-Cover with $n_1=n^{2/5}$. If we let $M$ be the resulting matching one can show that w.s.h.p. $D_2=\emptyset$ and $M$ saturates every vertex not in $V_2$. Therefter we can go straight to phase 2 of the algorithm to transform the paths induced by $M$ into a Hamilton cycle.

If $n^{-1/3}> p\geq \frac{101\log n}{n}$ then we greedily construct a path $P$ until $n'=\frac{30\log n}{p}$ vertices are left uncovered. We construct $P$ by first adding to it the vertex $v_1$. Then, while $|P|< n- \frac{30\log n}{p}$ we let $u$ be the last vertex added to $P$ and we query edges form $u$ to $[n]\setminus V(P$), until we identify an edge $e$ that belongs to $G$. Then, we extend $P$ using $e$. This process fails to construct a path of length $n-\frac{30\log n}{p}$ only if there exists a vertex for which $\frac{30\log n}{p}$ many queries were made none of which resulted to an edge. Thus it fails with probability at most $n\Pr(Bin(30\log n/p,p)=0)=o(n^{-7})$. Now, with $v_P(=v_1),v_P'$ being the endpoints of $P$, observe that $G'=G\setminus (V(P)\setminus \{v_P,v_P'\})$ is distributed as $G(n',p)$ with, $n'p\leq 30\log n\leq 90\log n'$ and so far we have queried in $G'$ only edges that are incident to $v_P$. We then run the algorithm described in this section with the twist that before running CoverGreedy we add the edge $\{v_P,v_P'\}$ to $M$. One can show that w.s.h.p. $|Dangerous|=\emptyset$ and with minor additional modifications one can get a Hamilton cycle of $G'$ that spans $\{v_P,v_P'\}$. Such a cycle corresponds to a Hamilton cycle of $G$.

If $p\leq p^*=\frac{0.99\log n}{n}$ then instead of running the algorithm found in this section we iteratively  query the entries of a row  for which no entry is equal to 1 until we reveal an entry that is equal to 1 or all the entries of  that row are revealed. If  we identify an all 0 row then we return ``G is not Hamiltonian". Else we return FAILURE. Hence, we read at most $n$ 1's and
$$\Pr(FAILURE) \leq (1-(1-p^*)^n)^{0.5n}\leq e^{-e^{-np^*} 0.5n}=o(n^{-7}).  $$
\bibliographystyle{plain}
\bibliography{bib}

\begin{appendix}
\section{Hamilton cycles in the randomly ordered adjacency lists model: Proofs of auxiliary lemmas
}\label{app:lists}

In this section we restate and prove Lemmas given in Section \ref{sec:adjlists}. Recall we let $G\sim G(n,p)$ for some $p=p(n)\geq 0$ and we assume that for $v\in [n]$ we are given a randomly ordered adjacency list $L_v$ consisting of the neighbors of $v$ in $G$. We have also defined $p^*=\frac{0.999\log n}{n}$, $q=\min\{0.25np, 100\log n\}$, $U_i=\{v\in [n]:(v,\text{ undefined})\in Q_{i}\}$, $D$ be set the vertices of degree at most $3.001q$ in $G$ and $B$ be the set of vertices that enters the set $V_7$ during Step 1 of Algorithm 1.

\begin{lemma}\label{lem:app:exposerandomness}
If $p\geq p^*$ then w.s.h.p. the following hold:
\begin{itemize}
    \item[(a)] $D\leq n^{0.97}$.
    \item[(b)]  For every $v\in[n]$ there exists at most 300 vertices in $D$ from $v$ within distance 10.
    \item[(c)] While at most $i< \min\{10^{5}n,\tau \}$ queries have been made, with $v=f(Q_{i})$, if $v\notin D$ then for $u\in [n]\setminus D$,
\begin{align*}
\Pr(q_{i}=(v,u)|Q_{i})
\begin{cases}
= (1+o(1)) \cdot \frac{1}{n} \cdot \frac{np-d_{Q_{i}}^+(u)}{np}  &\text{ if } (v,u),(u,v) \notin Q_{i}.
\\  \leq \frac{2000}{np} &\text{ if }  (u,v) \in Q_i
\\ = 0 &\text{ if } (v,u)\in Q_i.
\end{cases}
 \end{align*}
\item[(d)] $\tau > 10^5 n$.
\end{itemize}
\end{lemma}

\begin{proof}
(a)$|D| \geq n^{0.97}$ only if there exists a set of vertices $S$ of size $s=n^{0.97}$  such that every vertex in $S$ has at most $3.001q\leq 0.751np$ neighbors in  $ V\setminus S$. Therefore,
\begin{align*}
    \Pr(|D|\geq s) &\leq  \binom{n}{s} \bigg(\sum_{k\leq 3.01q} \binom{n}{k}(p)^{k}(1-p)^{n-k} \bigg)^s
\leq \bfrac{en}{s}^s \bigg( 4 \bfrac{enp}{0.751np}^{0.751np}
e^{-np+0.751np^2} \bigg)^s 
\\&= 
\exp\bigg\{s\bigg(O(1)+\log (n/s) + 1.287\cdot 0.751np-
 np+0.751np^2\bigg)\bigg\}  
=o(n^{-7}).
\end{align*}

(b) In the event that there exists  $v\in[n]$ and at least 301 vertices in $D$ within distance 10 from $v$, $G$ spans a tree $T$ on at most 3011 vertices such that there exists $S\subset V(T)$ of size 301 with the property that for $u\in S$, $u$ has at most $3.001q\leq  0.751np$ neighbors in $[n]\setminus S$. The probability of this event occurring is at most,
\begin{align*}
&\sum_{t=301}^{3011} \binom{n}{t} t^{t-2}p^{t-1} \binom{t}{301} \bigg(\sum_{k\leq 0.751np}\binom{n}{k}p^k(1-p)^{n-k-t}
\bigg)^{301}
\\&\leq 2^{3011}\sum_{t=301}^{3011} n(np)^{t-1} \bigg(4\bfrac{enp}{0.751np}^{0.751np}e^{-(1+o(1))pn} \bigg)^{301}
\\& = 2^{3011}n(np)^{3011} \exp\big\{301\cdot (1+o(1)) \cdot \big(1.287\cdot 0.751-1 \big) pn \big) 
\big\}
\\& \leq 
2^{3011}n\exp\big\{3011\log(np) -10np\big\}
=o(n^{-7}).
\end{align*}

(c) Let $Q=Q_i$ and $q_j=(a_j,b_j)$ for $j\leq i$. The conditions (C1), (C3) imposed on our algorithm imply that $d_{Q}^-(w)<\frac{\log n}{\log\log\log n}$ and $d_{Q}^+(w)\leq 3q$ for $w\in [n]$. For $x\in [n]$ let $\cG^x$  be the set of graphs on $[n]$ such that if $F\in \cG^x$ then (i) $F$ contains $G(Q)$ as a subgraph, (ii) $F$ spans at least $n-n^{0.97}-1$ vertices of degree at least $3.001q-1$ and (iii) for $z\in [n]$ if $(z,\text{undefined})\in Q$  or $z=x$  then the neighborhoods of $w$ in both $G(Q)$ and $F$ are identical. In addition let $\cG^x_s$ be the set of star-graphs on $[n]$ centered at $x$ that have at least $3.001q-d_{Q}(x)$ edges and are edge disjoint from $G(Q)$.

Let $0<\epsilon<0.1$ and $\mathcal{E}_x$ be the event that $G\setminus \{x\}$ spans at least $0.1n\log n$ edges. Part (a) implies that $\cE_x$ occurs w.s.h.p. We start by showing that conditioned on $Q$, $x\notin D$ and $\cE_x$ occurring, the degree of $x$ in $G$ lies w.s.h.p. in $I_p=[(1-\epsilon)np,(1+\epsilon)np]$. For $F\in \cG^x$ and $F_s\in \cG^x_s$ let 
$$ p(F)= p^{|E(F)|}(1-p)^{\binom{n-1}{2}-|E(F)|},  \hspace{5mm}      p(F_s)= p^{|E( F_s)|}(1-p)^{(n-1)-|E(F_s)|},$$
$$p_{F}= \prod_{j=1}^{i} \frac{\mathbb{I}(a_j\neq x)}{d_F(a_j)-d_{Q_{j-1}}^+(a_j)}, \hspace{5mm} p_{F,F_s,x}=\prod_{j=1}^{i} \frac{\mathbb{I}(a_j= x)}{d_{F\cup F_s}(a_j)-d_{Q_{j-1}}^+(a_j)},$$
$$p_{F,F_s}= \prod_{z\in N_{F_s}(x)} \frac{d_F(z)+1-d_{Q}^+(z)}{d_F(z)+1}$$
and lastly,
$$ w(F,F_s)= p(F)p(F_s) \cdot p_F\cdot p_{F,F_b} \cdot p_{F,F_s,x}.$$
Then,
\begin{align}\label{eq:Qi-1}
    \Pr(Q, x\notin D, \cE_x )&= \sum_{F\in \cG^x, F_s\in \cG^x_s} w(F,F_s). 
\end{align}

For $F\in \cG^x$ and $F_1,F_2\in \cG^x_s$ such that $F_2=F_1\cup \{x,z\}$ for some $z=z(F_2)\in [n]\setminus N_{G(Q)}(x)$ we have,
\begin{align}\label{eq:f1f2}
w(F,F_2)  = \frac{p}{1-p} \cdot \frac{d_{F\cup F_1}(x)+1-d_{Q}^+(x)}{d_{F\cup F_1}(x)+1 }\cdot \frac{d_F(z)+1-d_{Q}^+(z)}{d_F(z)+1} \cdot w(F,F_1).
\end{align}
$d^+_Q(x)\leq 3q$ and $d^-_{Q}(x)\leq \frac{\log n}{\log\log \log n}$ imply that $d_{F\cup F_1}(x)+1-d_{Q}^+(x)=(1+o(1))d_{F_1}(x)$. Thus  if $d_{F\cup F_1}(x) \geq (1+\epsilon/2)np$ then   \eqref{eq:f1f2} implies that
\begin{align*}
w(F,F_2)  \leq  \frac{p}{1-p} \cdot \frac{(1+o(1))d_{F_1}(x)}{d_{F_1}(w)+d^+_Q(x) } \cdot w(F,F_1).
\end{align*}
Therefore for
$(1+\epsilon/2)np- d_{Q}^+(x) \leq d$,  
\begin{align}\label{eq:Qidegree+}
\sum_{\substack{F\in \cG^x, F_2\in \cG^x_s:\\ d_{F_2}(x)=d+1}} w(F,F_2) &\leq \bfrac{(n- d_{Q}(x)-1)-d}{d+1} \bfrac{p}{1-p}  \frac{(1+o(1))d }{d+d_{Q}^+(v)} 
\sum_{\substack{F\in \cG^x, F_1\in \cG^x_s:\\d_{F_1}(x)=d}} w(F,F_1) \nonumber
\\& \leq  (1+o(1)) \bfrac{n-(1+o(1))d_{Q}^+(x)-d}{d+d_{Q}^+(x)} \bfrac{p}{1-p} \sum_{\substack{F\in \cG^x, F_1\in \cG^x_s:\\d_{F_1}(x)=d}} w(F,F_1) \nonumber
\\& \leq (1+o(1))\bigg(1-\frac{\epsilon/2}{(1+\epsilon/2)(1-p)}\bigg)  \sum_{\substack{F\in \cG^x, F_1\in \cG^x_s:\\d_{F_1}(x)=d}} w(F,F_1).
\end{align}

On the other hand each $F\in \cG_x$ spans at least $n-1-n^{0.97}$ vertices of degree at least $3.001q-1\geq 0.1\log n$. Moreover as $i<10^5n$ at most $n/\log\log n$ of those vertices have degree larger than $2\cdot 10^5\log\log n$ in $G(Q)$. Thus there exists a set $S$ of size $n-o(n)$ such that if $z\in S$ then $d_{Q}(z)=o( d_{F}(z) )$. In such a case \eqref{eq:f1f2} becomes,

\begin{align*}
w(F,F_2)  =  \frac{(1+o(1))p}{1-p} \cdot \frac{d_{F_1}(x)}{d_{F_1}(x)+d^+_{Q}(x)} \cdot w(F,F_1).
\end{align*}
Thus for $(1-\epsilon/2)np-d_{Q}^+(x) \geq d$, 

\begin{align}\label{eq:Qidegree-}
\sum_{\substack{F\in \cG^x, F_1\in \cG^x_s:\\d_{F_1}(x)=d}} w(F,F_1)
&\leq 
\bfrac{(1+o(1))d+1}{n-o(n)-d_{Q}(x)-d} \bfrac{(1-p)}{p}  \bfrac{d+d_{Q}^+(x)}{d}  \sum_{\substack{F\in \cG^x, F_2\in \cG^x_s:\\d_{F_2}(x)=d+1}}  w(F,F_2) \nonumber
\\& = \bfrac{(1+o(1))(d+d_{Q}^+(x))(1-p)}{[(1+o(1))n-d-d_{Q}^+(x)]p}  \sum_{\substack{F\in \cG^x, F_2\in \cG^x_s:\\d_{F_2}(x)=d+1}}  w(F,F_2) \nonumber
\\&\leq (1+o(1)) \bfrac{(1-\epsilon/2)np(1-p)}{[n-(1-\epsilon/2)np]p}   \sum_{\substack{F\in \cG^x, F_2\in \cG^x_s:\\d_{F_2}(x)=d+1}}  w(F,F_2) \nonumber
\\&\leq (1+o(1))\bigg(1- \frac{\epsilon/2}{1-(1-\epsilon/2)p}\bigg)   \sum_{\substack{F\in \cG^x, F_2\in \cG^x_s:\\d_{F_2}(x)=d+1}}  w(F,F_2). 
\end{align}
Now let $v=f(Q_{i})$  and let $u\in [n] \setminus D$. Equations \eqref{eq:Qi-1}, \eqref{eq:Qidegree+} and \eqref{eq:Qidegree-} imply that
\begin{align}\label{eq:Qi-1:restr}
    \Pr(Q_{i}, v\notin D,u\notin D, \cE_v,\cE_u)&=(1+o(1))\Pr(Q_{i},\cE_{v,u}),
    \end{align}
where $\cE_{v,u}$ is the event that $v,u \notin D$, $\cE_v,\cE_u$ occurs and $d_{G}(v), d_{G}(u)\in I_p=[(1-\epsilon)np, (1+\epsilon)np]$.
Now, assume that $(v,u),(u,v) \notin Q_{i}$. \eqref{eq:Qi-1:restr} gives,
\begin{align*}
    \Pr(Q_{i}, \cE_{v,u})&= (1+o(1))\sum_{\substack{F\in \cG^{v}, F_s \in \cG^{v}_s: (v,u)\notin E(F_s), \\   d_{F\cup F_s}(v), d_{F\cup F_s}(u) \in I_p
    }} w(F,F_s)\cdot \bigg(\frac{p  w(v,u,F,F_s) }{1-p}   +1 \bigg).
\end{align*}
where $w(v,u,F,F_s)=\bfrac{d_{F\cup F_s}(v)+1-d_{Q_i}^+(v)}{d_{F\cup F_s}(v)+1}\bfrac{d_{F\cup F_s}(u)+1-d_{Q_i}^+(u)}{d_{F\cup F_s}(u)+1}$. By considering the cases $p\leq n^{-0.5}$, $p\geq n^{-0.5}$ and using at the second case the inequalities  $d^+_{Q_i}(u), d^+_{Q_i}(v)\leq 3.001q \leq 400\log n$ and $d_{F\cup F_s}(v), d_{F\cup F_s}(u)\in I_p$, the above equation implies that
\begin{align}\label{eq:Qi-1:uv}
    \Pr(Q_{i},\cE_{v,u})&= (1+o(1))\sum_{\substack{F\in \cG^{v}, F_s \in \cG^{v}_s: (v,u)\notin E(F_s), \\   d_{F\cup F_s}(v), d_{F\cup F_s}(u) \in I_p
    }} w(F,F_s)\bigg( \mathbb{I}(p<n^{-0.5})+\frac{\mathbb{I}(p\geq n^{-0.5})}{1-p}\bigg).
\end{align}

Similarly,
\begin{align}\label{eq:Qi}
    \Pr(Q_{i},q_{i+1}=(v,u),\cE_{v,u})&= (1+o(1))\sum_{\substack{F\in \cG^{v}, F_s\in \cG^{v}_s:e\notin E(F_s), \\  d_{F\cup F_s}(v), d_{F\cup F_s}(u)\in I_p 
    }} w(F,F_s)\cdot w'(v,u,F,F_s)
\end{align}
where,
\begin{align*}
    w'(v,u,F,F_s) &= \frac{p}{1-p} \cdot  \frac{ w(v,u,F,F_s)}{(d_{F\cup F_s}(v)+1)-d_{Q_i}^+(v)}
=(1+o(1)) \frac{d_{F\cup F_s}(u) -d_{Q_i}^+(u)}{d_{F\cup F_s}(u) \cdot d_{F\cup F_s}(v)\cdot (1-p)/p}
\\&\in \bigg[\frac{(1+o(1))(1-\epsilon)np-d_{Q_i}^+(u)}{(1-\epsilon)^2n^2p(1-p)},\frac{(1+o(1)(1+\epsilon)np-d_{Q_i}^+(u)}{(1+\epsilon)^2n^2p(1-p)} \bigg].
\end{align*}
At the last line we used that $d_{F\cup F_s}(v), d_{F\cup F_s}(u)\in I_p$.
By taking $\epsilon \to 0$,
\eqref{eq:Qi-1:uv} and \eqref{eq:Qi} give, 
$$ P(q_{i+1}=(v,u)|Q_{i},\cE_{v,u})=
(1+o(1))\frac{np-d_{Q_i}^+(v)}{n^2p}.
 $$
Now if $(v,u)\in Q_i$ then the entry of $L_{v}$ corresponding to the edge $\{v,u\}$ has already been revealed  and $\Pr(q_{i+1}=(v,u)|Q_{i})=0$. 

Finally if $(v,u) \in Q_{i_i}$ then since $v\notin D$ at least $0.001np$ entries of $L_v$ have not been revealed. Thus $u$ is the next entry in $L_v$ with probability at most $\frac{1}{0.001np}\leq \frac{2000}{\log n}$.

(d) For $u\in [n]\setminus D$ and $i \leq \tau$, $d_{Q_i}(u)-d_{Q_i}^-(u)$ is at most equals to $|\{j\leq i: q_j=(v,u) \text{ and }(u,v)\notin Q_{j-1}\}|$ plus the number of neighbors of $u$ in $D$, which by part (b) is w.s.h.p. at most $300$.  As at most $0.75\log n$ entries from each list are queried, part (c) implies that the probability there exists $i\leq 10^{5}n$ and $u\in [n]\setminus D$ such that $d_{Q_i}(u)-d_{Q_i}^-(u)> \frac{\log n}{\log \log \log n}$ is at most 
\begin{align*}
    n\Pr\bigg(Bin\bigg(10^5 n, \frac{2 }{n}\bigg)\geq \frac{\log n}{\log \log \log n}-300 \bigg)
    &\leq n\binom{10^5n}{\frac{\log n}{\log \log \log n}-300} \bfrac{2 }{ n}^{\frac{\log n}{\log \log \log n}-300}
    \\&\leq n\bfrac{2\cdot e\cdot 10^5 }{\frac{\log n}{\log \log \log n}}^{\frac{\log n}{\log \log \log n}-300}=o(n^{-7}).
    \end{align*}
\end{proof}

\begin{lemma}\label{lem:app:phase1prep}
W.s.h.p. if $p^* \leq p \leq  \frac{1000\log n}{n} $ then the following hold:
\begin{itemize}
    \item[(a)] $|B\setminus D|\leq  n^{0.97}$.
    \item[(b)] For every $v\in[n]$ there exists at most 300  vertices in $B\setminus D$ within distance 10 from $v$ in $G_1$.
    \item[(c)] No subgraph $F$ of $G_1$ has at most $10^5$ vertices and spans more than $|V(F)|+8$ edges.  
    \item[(d)] For every $v\in V$ the number of vertices in $\cup_{i\in[4]}I_i$ within distance  at most 10 from $v$ is at most 300.
\end{itemize}
In the case $p\geq \frac{1000\log n}{n}$ we have that $B=D=\emptyset$.
\end{lemma}

\begin{proof}
(a) A vertex $v\in [n]$ belongs to $B\setminus D$ if there exists $j\in [4]$ such that out of the at least $0.25\log n$ queries made by Algorithm 1 less than 4000 belong to $V_j$.  Lemma \ref{lem:app:exposerandomness} implies that this occurs, independently of the rest of the queries, with probability at most
$$\Pr(Bin(0.25\log n,0.24)< 4000) \leq \binom{0.25\log n}{4000}(1-0.24)^{0.25\log n-4000}
\leq n^{-0.06}.
$$
Thus $\Pr(|B\setminus D|\geq n^{0.97})$ is bounded above by, 
\begin{align*}
\Pr(Bin(n,n^{-0.06})\geq n^{0.97}) \leq \binom{n}{n^{0.97}}(n^{-0.06})^{n^{0.97}} \leq \bfrac{en\cdot n^{-0.06}}{n^{0.97}}^{n^{-0.97}}=o(n^{-7}).
\end{align*}

(b) In the event that there exists  $v\in[n]$ with at least 301  vertices in $B\setminus D$ within distance 10 in $G_1$ then $G$ spans a tree $T$ on at most $t$ vertices such that there exists  $S\subset V(T)$ of size $301$   such that the following hold:
(i) $301\leq t\leq 3001$, (ii) with $Q_T=\{q_{i_1},q_{i_2},...,q_{i_{r}}\}$ being  the set of queries made during Algorithm 1 with their first coordinate in $V(T)$, hence $|Q_T|\leq qt=\min\{0.25npt,100t\log n\}$, at least $t-1$ element of $Q_T$ define an edge spanned by $V(T)$ and not defined by a previous query,
(iii) for every vertex $v\in S$ there exists $j\in [4]$ such that  $v$ has at most 4000 neighbors in $V_j\setminus S$ in $G_1$.
Lemma \ref{lem:app:exposerandomness} implies that the probability of the above occurring is bounded above by 
\begin{align*}
    &\sum_{t=301}^{3001} \binom{n}{t}\binom{t}{301}\binom{qt}{t-1} \bfrac{(1+o(1))t}{n}^{t-1} \bigg( \Pr(Bin(0.25\log n,24)<4000\bigg)^t
    \\& \leq \sum_{t=301}^{3001} \bfrac{en}{t}^t t^{301} (qt)^t\bfrac{t}{n}^{t-1}  \bigg( n^{-0.06}\bigg)^t =n^{1+o(1)-0.06\cdot 301}=o(n^{-7}).
\end{align*}

(c) In the event that $G_1$ spans a graph on $s\leq 10^5$ vertices and $s+9$ edges there exists  $S\in[n]$, $|S|=s$ such that with $Q_S=\{q_{i_1},q_{i_2},...,q_{i_{r}}\}$ being  the set of queries made during Algorithm 1 with their first coordinate in $S$, hence $|Q_S|\leq qs$, at least $s+9$ element of $Q_S$ define an edge spanned by $S$ and not defined by a previous query. Lemma \ref{lem:app:exposerandomness} implies that the probability of the above occurring is bounded above by 
\begin{align*}
    &\sum_{s=6}^{10^5} \binom{n}{s}\binom{qs}{s+9} \bfrac{(1+o(1))s}{n}^{s+9}=o(n^{-7}).
\end{align*}

(d) The probability that there exists  $v\in [n]$ and $j\in 4$ such that $v$ has more that 300 neighbors in $I_j$ in $G_1$, hence in $G$, is at most
$$4n\binom{20n^{0.97}}{301}p^{301} \leq 4n (n^{0.97}p)^{301} \leq 4n (n^{-0.03}\log n)^{301}
=o(n^{-7}).$$

Now let $p\geq \frac{1000\log n}{n}$. A vertex belongs to $D$ if it has degree less than $3.001q\leq 0.76np$. Thus,
$$\Pr(D\neq \emptyset )\leq n \sum_{i\leq 0.76np }\binom{n}{i}p^i(1-p)^{n-i}\leq n^2 \bfrac{enp}{0.76np}^{0.76np}e^{-0.99np}=o(n^{-7}).$$
Conditioned  on $D=\emptyset$ a vertex $v\in [n]$ belongs to $B$ if there exists $j\in [4]$ such that among the first $100\log n\leq q$ entries of $L_v$ at most $3999$ belong to $V_j$. Thus Lemma \ref{lem:exposerandomness} implies,
$$\Pr(B\neq \emptyset )\leq 4n \Pr(Bin(100\log n,0.24)\leq 3999) \leq 4n  (100\log n)^{3999} 0.76^{100\log n-3999}=o(n^{-7}).$$
\end{proof}

\begin{lemma}\label{lem:app:phase2prep}
W.s.h.p. the following hold:
\begin{itemize}
    \item[(a)] Every pair of sets $X'\subset X, Y'\subset Y$ of size at least $0.002n$ spans an edge in $H_1$.
    \item[(b)] For every $S\subset X$ or  $S\subset Y$, $|S| \geq 10^{-6}n$ implies that the neighborhood of $S$ in $H_1$ has size at least  $ 0.002n$.
    \item[(c)] For every $S\subset X$ (resp. $S\subset Y$) of size at most $10^{-6}n$ the set $\{u: u\in L_v^1 \text{ for some }v \in S\}$ has size at least $2000|S|$.
\end{itemize}
\end{lemma}
\begin{proof}
(a)  The probability that there exists $X'\subset X$ of size $n_1=0.002n$ such that $\{v\in L_x^1:x\in X''\}$ has size at most $|Y|-0.002n$ is bounded above by,
\begin{align*}
    &\binom{|X|}{n_1}\binom{|Y|}{n_1} \bfrac{\binom{|Y|-n_1}{3000}}{\binom{|Y|}{3000}}^{n_1}
    \leq \bigg( (1+o(1))\bfrac{0.25en}{0.002n}\bfrac{0.248n}{0.25n}^{1500} \bigg)^{0.004n} =o(n^{-7}).
\end{align*}

(b) In the event that there exists $W\in \{X,Y\}$ and $W'\subset W$ such that $|W'|= 10^{-6}n$ and $|N_{H_1}(W')| \leq n_1$  the set $\{w\in L_w^1:w\in W'\}$ has size at most $n_1$. With $z=|X|=|Y|=(1+o(1))0.25n$, the probability of this event occurring is bounded above by,
\begin{align*}
    2\binom{z}{10^{-6}n}\binom{z}{n_1}\bfrac{\binom{n_1}{3000}}{\binom{z}{3000}}^{10^{-6}n}
    \leq 2\bigg( (0.25e10^6+o(1))(125)^{2000}\bfrac{1}{125}^{3000} \bigg)^{10^{-6}n}=o(n^{-7}).
\end{align*}

(c) With $z=|X|=|Y|=(1+o(1))0.25n$,
\begin{align*}
    \Pr(c) &\leq 2\sum_{0\leq i\leq 10^{-6}n}\binom{z}{i}\binom{z}{2000i}\bfrac{\binom{2000i}{3000}}{\binom{z}{3000}}^{z}
    \leq 2\sum_{0\leq i\leq 10^{-6}n} \bigg(\frac{ez}{i}\bfrac{ez}{2000i}^{2000} \bfrac{2000i}{z}^{3000}  \bigg)^i
    \\&\leq 2\sum_{0\leq i\leq 10^{-6}n} \bigg( 2000e \cdot \bfrac{2000ei}{z}^{999} \bigg)^i 
    =o(n^{-7}).
\end{align*}

\end{proof}

\begin{lemma}\label{lem:app:endpoints}
W.s.h.p. there exists sets $S(v)\subset N_{G_1}(v)$ for $v\in V_R$ such that $|S(v)|\geq 10^4$  and no vertices in $V_R\cup \big(\cup_{v\in R} S(v)\big) $ are within distance 10 in $H_0$.
\end{lemma} 

\begin{proof}
For each vertex $v\in V_R$ has at least 16000 neighbors in $G_1$, 4000 in each $V_j'$, $j\in [4]$. Part (c) of Lemma \ref{lem:phase1prep} implies that there further exists a set $N(v)$ consisting of 15000 of these neighbors such that no 2 vertices in $N(v)$ lie within distance 1000 (defined by $H_0$) or are within $H_0$-distance 500 from $v$. One can generate the sets in question by starting with $S=V_R$ and iteratively finding and removing from $S$ a vertex $v$ which has a set of $10^4$ neighbors $N$ each at distance (defined by $H_0$) at least 10 from $\cup_{u\in S}N(u)$ and setting $S(v)=N$.

Consider the auxiliary graph $G'$ on $S \cup  \big( \bigcup_{s\in S}N(s)\big)$ where $u\in S$ is connected to every vertex in $S(u)$ and for $u,v\in S$, $w\in S(v)$ and $z\in S(u)$, $wz\in E(G_1)$ if $d_{H_0}(w,z)\leq 10$. Let $\cE$ be the event that the above procedure does not outputs a set $S(v)$ for every $v\in V_R$. If $\cE$ occurs then there exists 
$S'\subset S$ such that for $s\in S'$ at least $5001$ vertices in $N(S)$ have at least 2 neighbors in $G'$. Hence $G'$ spans a subgraph $G''$ with $15001s$ vertices and $17500s$ edges for some $1\leq s\leq \log^2 n$. Due the construction of $G'$ this corresponds to a subgraph $G_1'$ of $G_1$ that spans $15001s+k$ vertices and $17500s+k$ edges for some $k\leq 10\cdot 17500s$.  
Thus if $p^*\leq p \leq \frac{\log^9n}{n}$ then the probability that $G_1$, hence $G$ spans such a subgraph is bounded by
\begin{align*}
    \Pr(\cE)&\leq \sum_{s=1}^{\log^2 n} \sum_{k=0}^{175000s} \binom{n}{s+k} \binom{\binom{15001s+k}{2}}{17500s+k}p^{17500s+k}=o(n^{-8}).
\end{align*}
On the other hand, by Chernoff Bound,  if $p \geq \frac{\log^9n}{n}$ then w.s.h.p. every vertex has degree at least $0.5np$. If $\cE$ occurs then each edge $\{u,v\}$ of $G_1'$ satisfies either $u\in \{L_v(i):i\leq 3q\}$ or $v\in \{L_u(i):i\leq 3q\}$. Therefore, as $q\leq 100\log n$ 
\begin{align*}
    \Pr(\cE)&\leq \sum_{s=1}^{\log^2 n} \sum_{k=0}^{175000s} \binom{n}{s+k} \binom{\binom{15001s+k}{2}}{17500s+k}p^{17500s+k}\bfrac{3q}{0.25np}^{17500s+k}=o(n^{-8}).
\end{align*}

\end{proof}

\section{Matrices}\label{app:matrix}
The Lemma below refers to the set $Dangerous$ that we get after the execution of SatDangerous.
\begin{lemma}\label{lem:app:satdanger:basicprop}
W.s.h.p. the following hold:
\begin{itemize}
\item [(i)] $\Delta(G)\leq 200\log n$.
\item[(ii)] There does not exist a connected subgraph of $G$ spans $k$ vertices and more than $k+8$ edges for $k\leq 10^5$.
\item[(iii)] $|Dangerous|\leq n^{0.1}$.
\item[(iv)] No connected subgraph of $G$ spans $k\leq 1000$ vertices out of which more than 9 vertices belong to $Dangerous$.
\end{itemize}
\end{lemma}
\begin{proof}
(i) $\Pr(\Delta(G)> 200\log n)\leq n\Pr(Bin(n,(100\log n)/n)>200\log n)=o(n^{-7})$.

(ii) \begin{align*}
    \Pr((ii)) \leq \sum_{k=1}^{10^5}\binom{n}{k} \binom{k^2}{k+9} \bfrac{100\log n}{n}^{k+9}=o(n^{-7}).
\end{align*}

(iii) If $v\in Dangerous$ then during one of the algorithms GreedyCover and SatDangerous w.s.h.p. we queried a set of at least $0.99n$ edges incident to $v$ that had not been queried earlier and at most $0.01\log n$ of those belong to $G$. Therefore,
\begin{align*}
\Pr(|Dangerous|\geq n^{0.1})& \leq \binom{n}{n^{0.1}}\bigg(\Pr(Bin(0.99n,p)<0.01n)\bigg)^{n^{0.1}}
\\&\leq \bigg(en^{0.1}\Pr(Bin(0.99n,p)<0.01n)\bigg)^{n^{0.1}}=o(n^{-7}).
\end{align*}

(iv) Let $\cB$ be the event that $G$ spans a connected subgraph on $k\leq 1000$ vertices out of which more than 9 vertices belong to $Dangerous$ then it spans a tree on $k$ vertices. If $\cB$ occurs then $G$ spans a tree $T$ such that $V(T)=k\leq 1000$ and $V(T)$ contains a set $B$ of 10 vertices with the following property. For $v\in B$ during one of the algorithms GreedyCover and SatDangerous w.s.h.p. we queried a set of at least $0.99n$ edges incident to $v$, not spanned by $V(T)$, that had not been queried earlier and at most $0.01\log n-|B|$ of those belong to $G$. Thus,
\begin{align*}
\Pr(\cB)\leq \sum_{k=1}^{1000} \binom{n}{k}k^{k-2}p^{k-1} \binom{k}{10}\bigg(\Pr(Bin(0.99n,p)<0.01\log n-k)\bigg)^{10}=o(n^{-7}).
\end{align*}
\end{proof}

\section{The Karp-Sipser algorithm}\label{sec:app:KS}
In this section we prove Lemmas \ref{lem:special:gnm2} and \ref{lem:karp} which we now recall. 
\begin{lemma}\label{lem:app:special:gnm2}
W.s.h.p. $|Special|\leq n^{0.3}$.
\end{lemma}
\begin{lemma}\label{lem:app:karp}
W.s.h.p. $|Special_0|\leq n^{0.41}/3$.
\end{lemma}
For that observe that the subgraph of $G$ spanned by $V'$ as it is defined before the execution and upon the termination of SatDangerous is distributed according to $G_{|V'|,|E(G[V'])|}^{\delta\geq 2}$. For integers $n,m$ with $m\geq n$ we let $G_{n,m}^{\delta \geq 2}$ be a random graph that is  chosen uniformly from all the graphs on $[n]$ with $m$ edges and minimum degree 2. What follows is taken (with minor modifications) from \cite{anastos2021k}, \cite{anastos2020}, \cite{Hamd3}.

\subsection{Random Sequence Model}\label{refined}
We now take some time to explain the model that we use to generate and analyse $G_{N,M}^{\delta \geq 2}$. We use a variation of the pseudo-graph model of Bollob\'as and Frieze \cite{Bollfr} and Chv\'atal \cite{Ch} to generate $G_{N,M}^{\delta \geq 2}$. Given a sequence $\bx = (x_1,x_2,\ldots,x_{2M})\in [N]^{2M}$ of $2M$ integers between 1 and $N$ we can define a (multi)-graph
$G_{\bx}=G_\bx(N,M)$ with vertex set $[N]$ and edge set $\{(x_{2i-1},x_{2i}):1\leq i\leq M\}$. The degree $d_\bx(v)$ of $v\in [N]$ is given by 
$$d_\bx(v)=|\set{j\in [2M]:x_j=v}|.$$
If $\bx$ is chosen randomly from $[N]^{2M}$ then $G_{\bx}$ is close in distribution to $G_{N,M}$. Indeed,
conditional on being simple, $G_{\bx}$ is distributed as $G_{N,M}$. To see this, note that if $G_{\bx}$ is simple then it has vertex set $[N]$ and $M$ edges. Also, there are $M!2^M$ distinct equally likely values of $\bx$ which yield the same graph. 

We will use the above variation of the pseudo-graph model to analyze Karp-Sipser and the part of SatTrouble that mimics the Karp-Sipser algorithm.  As Karp-Sipser progresses vertices become matched, edges are deleted and vertices of small degree are identified.  As such we will need to impose additional constrains on the vertex degrees and our situation becomes more complicated. At any step of the algorithm we keep track of 2 sets $J_2$ and $J_1$ that partition the current vertex set, say $[N]$. $J_2$ is a set of vertices of degree at least 2 and it consists of vertices that have not been matched yet.  $J_1$ consists of the remaining vertices i.e. of vertices of degree 1 that have not been matched and it will be proven to be of size  $o(N)$.

So we let
$$[N]^{2M}_{J_2}=\{\bx\in [N]^{2M}:d_\bx(j)\geq 2\text{ for }j\in J_2, \text{ and } d_\bx(j)=1 \text{ for } j\in [N]\setminus J_2\}.$$
Let $G=G(N,M,J_2)$ be the multi-graph $G_\bx$ for $\bx$ chosen uniformly from $[N]^{2M}_{J_2}$. 
What we need now is a procedure that generates $G_\bx$ conditioned on $G_\bx$ being simple or equivalently a way to access the degree sequence of elements in $[N]^{2M}_{J_2}$. Such a procedure is given in \cite{Hamd3} and it is justified by Lemmas \ref{lem3} and \ref{lem:degDistributions} that
follow. Lemma \ref{lem3} states that the degree sequence of a random element of 
$[N]^{2M}_{J_2}$ (restricted to the set $J_2$) has the same distribution as the joint distribution of $\cP_1,\cP_2,...,\cP_{|J_2|}$
where (i) for $i\in J_2$, $\cP_i$ is a $Poisson(\lambda)$ random variable condition on being at least $2$  for some carefully chosen value of $\lambda$ and (ii) $\sum_{i=1}^{|J_2|}\cP_i=2M-(N-|J_2|)$. This fact is used in Lemma \ref{lem:degDistributions}
which establishes concentration of the number of vertices of degree $k$ in $J_2$. For the proofs of Lemmas \ref{lem3} and \ref{lem:degDistributions} see \cite{anastos2020}, \cite{anastos2021k}, \cite{Hamd3}.
We let
$$f_k(\l)=e^\l-\sum_{i=0}^{k-1}\frac{\l^i}{i!}$$
for $k\geq 0$.
\begin{lemma}
\label{lem3}
Let $\bx$ be chosen randomly from $[N]^{2M}_{J_2}$. Let $Z_j\,(j\in [J_2])$ be independent copies of a {\em truncated Poisson} random variable $\cP_i$, where
$$\Pr(\cP_i=t)=\frac{{\l}^t}{t!f_i({\l})},\hspace{1in}t=2,3,\ldots\ .$$
Here ${\l}$ satisfies
\begin{equation}\label{21}
\frac{{\l}f_{1}({\l})}{f_2({\l})}|J_2|=2M-(N-|J_2|).
\end{equation}
For $j\in J_0=[N]\setminus J_2$ let  $Z_j=1$. Then $\{d_\bx(j)\}_{j\in [N]}$ is distributed as $\{Z_j\}_{j\in [N]}$ conditional on $Z=\sum_{j\in [n]}Z_j=2M$.
\end{lemma}

\begin{lemma}\label{lem:degDistributions}
Let $\bx$ be chosen uniformly from $[N]^{2M}_{J_2}$ and $Z$ be a  truncated Poisson($\l$)  where
${\l}$ satisfies \eqref{21}.
 For $j\geq 2$ let $ \nu_{j}(\bx)$ be the number of elements $v\in J_2$ such that $d_\bx(v)=j$.
Then with probability $1-o(N^{-25})$,
$$\bigg|\nu_{j}(\bx)- |J_2|\Pr(Z=j)\bigg| \leq N^{1/2}\log^2 N  \text{ for } 2 \leq j\leq \log^2 N.$$
\end{lemma}

Now let $N,M\in \mathbb{N}$ with $N\leq M\leq N \log^{0.45}N$, $J_2=[N]$ and $\lambda$ be defined by \eqref{21}. Let $\cE$ be an occupancy event in $G_{N,M}^{\delta \geq 2}$. Denote by $G_{N,M}^{\delta \geq 2, seq}$ the random graph that is generated from the random sequence model by first choosing a random element of $[N]_{[N]}^{2M}$ and then generating the corresponding graph. Also denote by $G_{N,M}^{\delta \geq 2,Po(\lambda)}$ the random graph that is generated by first generating $N$ independent $Poisson(\lambda)$
random variables, conditioned on being at least 2, $P_1, P_2, ..., P_N$, then choosing a random sequence in $[N]^{\sum_{i\in[N]}P_i}$ with degree sequence $P_1, P_2, ..., P_N$ and finally generating the corresponding graph if $\sum_{i\in [N]} P_i$ is even. Then,
\begin{align}\label{model1}
\Pr\big(G_{N,M}^{\delta\geq 2}\in \cE\big) &\leq \Pr^{-1}\big(G_{N,M}^{\delta\geq 2,seq} \text{ is simple }\big) \Pr\big(G_{N,M}^{\delta\geq 2,seq}\in \cE\big)
\end{align}
and
\begin{align}\label{model2}
\Pr\big(G_{N,M}^{\delta\geq 2,seq}\in \cE\big) =
\Pr\bigg(G_{N,M}^{\delta\geq 2,Po(\lambda)}\in \cE \big| \sum_{i\in [N]} P_i=2M\bigg)
\leq O(N^{0.5})\Pr\big(G_{N,M}^{\delta\geq 2,Po(\lambda)}\in \cE \big).
\end{align}
For the equality in \eqref{model2} we used Lemma \ref{lem3} which implies that the models $G_{N,M}^{\delta\geq 2,seq}$ and  $G_{N,M}^{\delta\geq 2,Po(\lambda)}$ conditioned on $\sum_{i\in [N]} P_i=2M$ define the same distribution. 
\vspace{5mm}
\\We let $V''$ be equal to the set $V'$ exactly before the execution of SatDangerous. Also we let $n'=|V''|$ and $m'=|E(G[U])|$. Thus $G[V'']\sim G_{n',m'}^{\delta\geq 2}$. Let $G''\sim G_{n',m'}^{\delta \geq 2, seq}$.
\begin{lemma}\label{lem:g1'}
With probability $1-o(n^{-10})$, $\log^{0.35}n\leq \frac{m_1}{n_1} \leq \log^{0.44}n$ and $\Delta(G'')\leq \log n$. In addition,
$$  \Pr(G''\text{ is simple })\geq n^{-0.01}.$$
\end{lemma}
\begin{proof}
Lemma \ref{lem:greedycover} states that $5n_1\leq |U|\leq |V_2| \leq 10n_1$ and $|U|^2p=\Theta(|U|\log^{0.4}n)$. The Chernoff bound implies,
\begin{align}\label{eq:critical}
\Pr\bigg(\frac{m_1}{n_1} \notin (\log^{0.39}n,\log^{0.44}n) \bigg) \leq \binom{n}{n_1}e^{-\Theta\bfrac{n}{\log^{0.2} n }}=o(n^{-10}).
\end{align}
Now let $\lambda$ be such that $\frac{\lambda f_1(\lambda)}{f_2(\lambda)}=\frac{m'}{n'}$. Hence w.s.h.p. $\lambda \leq \log^{0.45} n$. Thereafter,\eqref{model2} implies,
$$\Pr(\Delta(G'')\geq \log n)\leq n^2\Pr(Poisson(\log^{0.45} n) \geq \log n) \leq \frac{2e^{-2\log^{0.45} n} (2\log^{0.45} n)^{\log n}}{ (\log n)!}= o(n^{-10}).$$
To prove that $\Pr(G''\text{ is simple })\geq n^{-0.01}$ we used a result of McKay \cite{mckay}. It implies that with $M_2=\sum_{i\in V''} d_\bx(i)(d_\bx(i)-1)$ and $\rho=M_2/m_1$, $\Pr(G'' \text{ is simple })=e^{-  (1+o(1)) \rho(\rho+1)}$. In our case one can use Lemma \ref{lem:degDistributions} to show that $\rho \leq 0.9\log^{0.45}n$ with probability at least $0.5n^{-0.01}$. Thus,
$$\Pr(G'' \text{ is simple })\geq 0.5n^{-0.01}+e^{- \log^{0.9}n}\geq n^{-0.01}.$$
\end{proof}

\paragraph{Proof of Lemma \ref{lem:app:special:gnm2}:}
Lemma \ref{lem:satdangerous} states that w.s.h.p. $|Special'|\leq 6n^{0.1}$.
In the event $|Special|>n^{0.3}$ and $|Special'|\leq 6n^{0.1}$ there exists $W'\subset Special$ such that (i) $|W'|=n^{0.3}$, (ii) $Special'\subset W$, (iii) every vertex in $W'\setminus Special'$ has at most 1 neighbor in $U$ that does not belong to its $G[W']$-component  and (iv) each component of $G[W']$ contains at least 1 vertex from $Special'$. Thus, there exists disjoint sets $S,T,Z\subset U$ such that $|S|\leq 6n^{0.1}$, $|T|\geq n^{0.2}/6$, $|Z|\leq |T|$, $G[S\cup T \cup Z]$ is connected and there is no edge from $T$ to $U\setminus (S\cup T\cup Z)$.

Now given $S, T, Z$ as above and a fixed tree $T'$ spanned by $S\cup T\cup Z$ the probability that $E(T')\subset E(G'')$ and $N_{G''}(T)\subset S\cup T\cup Z$ conditioned on with $\Delta(G'')\leq \log n$ is bounded by 
\begin{align}\label{tree}
&\prod_{uv\in E(T')}\bfrac{d(u)d(v)}{2|E(G')|}\bigg(1-\frac{2|E(G'')|-\sum_{v\in S\cup T\cup Z}d(v)}{2|E(G'')|} \bigg)^{\sum_{v\in T}d(v)-2|E(T')|} \nonumber
    \\&\leq \bfrac{(\log n)^2}{2n_1}^{|S|+|T|+|Z|-1}\bigg(1-\frac{2n_1-\log n(|S|+|T|+|Z|)}{2n_1} \bigg)^{\sum_{v\in T}d(v)-6|T|}+n^{-10}.
\end{align}
At the inequality we used Lemma \ref{lem:g1'}. On the other hand the probability that the sum of $|T|$ independent random variables $Poisson(\lambda)$ is smaller than $2|T|\log^{0.39} n$ is $o(n^{-9})$ for $|T|\geq n^{0.2}/6-6n^{0.1}$. Therefore equations \eqref{model1},\eqref{model2} and \eqref{tree} together with Lemma \ref{lem:g1'} imply,
\begin{align*}
    &\Pr(|Special|>n^{0.3}) \leq n^{0.6} \bigg\{ \sum_{s=1}^{6n^{0.1}}\sum_{t=n^{0.2}/6-s}^{n^{0.3}}\sum_{z=1}^{t} \binom{n}{s+t+z}\binom{s+t+z}{s}\binom{t+z}{t} (t+z+s)^{t+r+z-2} 
    \\& \times \bfrac{\log^2 n}{2n_1}^{s+t+z-1} \bfrac{(s+t+z)\log n}{2n_1} ^{1.5t\log^{0.39} n}+o(n^{-9})\bigg\}
    \\& \leq n^{0.6}  \sum_{s=1}^{6n^{0.1}}\sum_{t=n^{0.2}/6-s}^{n^{0.3}}\sum_{z=1}^{t} n \bfrac{en\log^2 n}{n_1}^{s+t+z}\bfrac{3et}{s}^s 2^{2t} \bfrac{3t\log n}{2n_1} ^{1.5t\log^{0.39} n}+o(n^{-8})
    \\& \leq 6n \sum_{t=n^{0.2}/7}^{n^{0.3}} n \bfrac{en\log^2 n}{n_1}^{3t}\bfrac{3et}{s}^s 2^{2t} \bfrac{3t\log n}{2n_1} ^{1.5t\log^{0.39} n}+o(n^{-8})=o(n^{-7}).
\end{align*}
\qed

Let $\tau$ be the number of times the while-loop at line 18 of SatTrouble or at line 1 of Karp-Sipser is executed. For $1\leq i\leq \tau$ by at time $i$ we refer to the moment exactly before the $i$th execution of the while-loop at line 18 of SatTrouble or at line 1 of Karp-Sipser. For $1\leq i\leq \tau$ we let $G_i=G[V_i]$, where $V_i$ equals to $V'$  at time $i$. We now let,
\begin{itemize}
\item $n_i$: size of $V_i$.
\item $m_i$: size of $E(G_i)$. \item $Z^i$: the set $V_i\setminus C_1$.
\item $Z_1^i$: the set $C_1$ at time $i$.
\item $Z_j^i$: the set of vertices in $V_i$ of degree $j$ in $G_i$, $J\geq 2$.
\item $\zeta_i=|Z^i_1|$.
\item $\l_i$ is the unique constant satisfying
\begin{equation*}
\frac{{\l_i}f_{1}({\l_i})}{f_2({\l_i})}|Z|=2m_i-(n_i-|Z_1^i|).
\end{equation*}
\end{itemize}
The following Lemma allows us to treat $G_i$ as a random graph of minimum degree 2.  For its proof see \cite{anastos2021k}, \cite{anastos2020}, \cite{Hamd3}.
\begin{lemma}\label{lem:markov}
If $G_1\sim G_{n_1,m_1}^{\delta\geq 2,seq}$ then $G_i\sim G_{n_i,m_i}^{\delta\geq 2,seq}$ for $1\leq i \leq \tau$. 
\end{lemma}
Let $\epsilon=10^{-9}$.
\begin{lemma}\label{lem:zeta}
With probability $1-o(n^{-8})$, if $G_1\sim G_{n_1,m_1}^{\delta\geq 2,seq}$ then for $i<\tau$ we have $\zeta_i< n^{0.4+\e}$.
\end{lemma}
We proof Lemma \ref{lem:zeta} in the next subsection. We now proceed with the poof of Lemma \ref{lem:app:karp}.
\paragraph{Proof of Lemma \ref{lem:app:karp}.}
Recall $G_1\sim G_{n_1,m_1}^{\delta \geq 2}$. Let $G_1'\sim G_{n_1,m_1}^{\delta \geq 2}$ be a random multigraph on $V'$  and $SSpecial_0$ be the set $Special_0$ that we get upon substituting $G_1$ by $G_1'$ at the beginning of SatTrouble. 

For every vertex $v\in SSpecial_0$  there exists a step $i$ such that either (i) $v\in Z_1^i$ and at step $i$ a neighbor of $v$ is matched and then removed from $V_i$ or (ii) $v\notin Z_1^i$, $1$ or 2 neighbors of $v$ are matched and then removed from $V_i$ and as a result $d_{G_i}(v)\geq 2$ edges incident to $v$ are removed. If the above occurs then we say that step $i$ witnesses an increase of $|SSpecial_0|$. 
Thus  step $i$ witnesses an increase of $|SSpecial_0|$ with probability at most   $2\Delta(G_1')\zeta_i/2m_i + O(1/m_i)$; the corresponding increase is at most $2\Delta(G_1')$.  If $|Special_0|$ reaches $n^{0.41}/3$ before time $\tau$ then in the event $\Delta(G_1')\leq \log n$ there are exists an integer $1\leq r\leq 1/\e^2$ and at least $\e^2 n^{0.4+2\e}$ steps with $m_i\in [n^{0.4+2\e+(r-1)\epsilon^2}, n^{0.4+2\e+r\epsilon^2}]$  that witness an increase of $|Special_0|$. The probability that this occurs, conditioned on $\Delta(G_1')\leq \log n$,  for a fixed $r$, while $\z_i\leq n^{0.4+\e}$,  is bounded by
\begin{align*}
\binom{n^{0.4+2\e+r\e^2}}{ \e^2 n^{0.4+2\e} } \bfrac{2n^{0.4+\e}\log n}{n^{0.4+2\e+(r-1)\e^2}}^{ \e^2 n^{0.4+2\e} } & \leq \brac{\frac{en^{r\e^2}}{\e^2 }  \cdot \frac{2\log n}{n^{\e+(r-1)\e^2}} }^{\e^2 n^{0.4+2\e}} = o(n^{-10}).
\end{align*} 
Lemma \ref{lem:g1'} gives that $\Pr(\Delta(G_1')> \log n)=o(n^{-10}).$
Thus \eqref{model1} implies,
\begin{align*}
    \Pr(|Special_0|\geq n^{0.41}/3)\leq n^{0.01}\cdot n\cdot o(n^{-10})=o(n^{-7}).
\end{align*}
\qed

\subsection{Proof of Lemma \ref{lem:zeta}}\label{subs:lem} For $1\leq i\leq \tau$ we let $\cH_i=\{(n_j,m_j)\}_{j\leq i}$ and assume that $G_1\sim G_{n_1,m_1}^{\delta\geq 2,seq}$. Lemma \ref{lem:markov} implies that $G_i\sim G_{n_i,m_i}^{\delta\geq 2,seq}$ for $1\leq i\leq \tau$.
We study $\z_i$ in three regimes (captured by the events  $\cA_i,\cB_i,\cC_i$ defined below) depending on whether we can approximate $Z_j^i, j\geq 2$ using the corresponding truncated Poisson distribution described earlier. In all the regimes, to bound $\z_i$ we show that it can be scholastically dominated by either a random walk with significantly large negative drift or a lazy random walk whose drift is asymptotically 0. The first kind of behavior is associated with the study of the random variables $X_i,Y_i'$ defined below while the second one is associated with $W_i$.

For $i<\tau$, we define the events
\[ 
\cA_i=\set{\lambda_i\geq m_i^{-0.2}}, \hspace{5mm} \cB_i=\set{m_i^{-0.2}>\lambda_i\geq |Z^i|\log^4 n}\]
\[\text{ and }\cC_i=\set{\lambda_j<  |Z^j|\log^4 n \text{ for some }j\leq i}.
\]
Also we let $\sigma=\min\{i:\cC_i \text{ occurs}\}$. For $i<\tau$, we also define the following random variables:
\begin{align*}
X_i&=(\zeta_{i+1}-\zeta_i) \mathbb{I}(\cA_i, 0<\zeta_i<n^{0.4+\e}).\\
 Y_i&=(\zeta_{i+1}-\zeta_i) \mathbb{I}(\cB_i, 0<\zeta_i<n^{0.4+\e}).
\\W_i&=(\zeta_{i+1}-\zeta_i) \mathbb{I}(i\leq \sigma, 0<\zeta_i<n^{0.4+\e})
\end{align*}
We will prove that for $0<i<\tau$ with probability $1-o(n^{-9})$
\begin{align}\label{goaly}
\min\{ \zeta_{i}, n^{0.4+\e}\} \leq M+\sum_{j=1}^{i-1}(X_i+Y_i+W_i)
\end{align}
where $M=\log^2 n$ is such that the following holds: with probability $1-o(n^{-9})$ for every $i\geq 0$ with $\zeta_i=0$ we have that $\zeta_{i+1} \leq M$. Our bound for $M$ is justified by the fact that the maximum degree in $G_i$ is at most $\log n$ with probability $1-o(n^{-10})$ (this follows from Lemma \ref{lem:g1'}). Similarly,  we get
\begin{equation}\label{eq:max}
\max\{X_i,Y_i,W_i:1\leq i\leq \tau\} \leq \Delta(G_1)\leq \log n \text{ with probability }1-o(n^{-10}).
\end{equation}

We use the inequality $i<\tau$, hence $m_i\geq 2n_{\tau}\geq 2n^{0.41}$, to impose that if $\zeta_i\leq n^{0.4+\e}$ then almost all of the vertices belong to $Z^i$. We will see from the analysis below that with probability $1-o(n^{-10})$,
\begin{equation}\label{important}
m_i\geq 2n^{0.41}\text{ implies }\zeta_i\leq n^{0.4+\e}.
\end{equation}

Equation (80) of \cite{Hamd3} states that, assuming one of the events $\cA_i,\cB_i$ occurs, 
\begin{equation}\label{neg}
\zeta_i>0\text{ implies }\E(\zeta_{i+1}-\zeta_i\mid \cH_i)\leq -\Omega(\min\set{1,\lambda_i}^2)+ O\bfrac{\log^2m_i}{\lambda_im_i}.
\end{equation}
In the following cases we will assume that $i<\tau$ and $\zeta_i>0$. The case $\z_i=0$ is handled by $M$ of \eqref{goaly}.

\textbf{Case 1: $\cA_i$ occurs.}
\\ If $\lambda_i\geq m_i^{-0.2}$ then we have from \eqref{neg} that  
\[
\E(X_i|\cH_i) \leq - c\lambda_i^2 \leq -c_1n^{-0.4}
\]
for some constant $c_1>0$.
Thus Azuma inequality gives,
\[
\sum_{j\geq 1} \Pr\bigg( \sum_{i=1}^j X_i \geq n^{0.4+\e/2} \bigg)
 \leq m_1 \max_{0\leq j\leq m_1}  \exp \bigg\{- \frac{(n^{0.4+\e/2}+ c_1jn^{-0.4})^2}{2j\log^2n}  \bigg\} +n^{-10} =o(n^{-10}).
\]
The $n^{-10}$ term accounts for using \eqref{eq:max} to bound $\max X_i$.

\textbf{Case 2: $\cB_i$ occurs.}
\\ To bound $\sum_{i\geq 1}^j Y_i$, let $R_i$ be the indicator of the event that $\zeta_i\neq \zeta_{i+1}$. Hence $Y_i=Y_iR_i$.

In the event $\cB_i$, $n_i\log^4 n \leq \lambda_i \leq m_i^{-0.2}$. In this case for $i\geq 2$, $|Z_i|$ is approximately equal to the sum of $|Z_i|$ independent random variables that follow Poisson($\lambda_i$) conditioned on having value at least 2. More precisely, it follows from Lemma 3.3 of \cite{Hamd3} that  as long as $\cB_i$ holds, we have 
\begin{equation}\label{1}
  \begin{split}
  \frac{|Z_3|}{|Z_2|}&= \frac{\lambda_i}{3}\brac{1+O(m_i^{1/2}\lambda_i\log^2m_i)},\\
 \frac{|Z_4|}{|Z_2|}&= \frac{\lambda_i^2}{12}\brac{1+O(m_i^{1/2}\lambda_i\log^2m_i)},\\
\sum_{i\geq 5} |Z_i|& \leq |Z_2| \lambda_i^3.
    \end{split}
\end{equation}
Recall that if $\zeta_i>0$ then the algorithm will  choose a vertex  $v\in Z_1^i$ and it will match it to some vertex $w$. Thus initially $\zeta_{i}$ will decrease by 1. 

For $w\in Z$ let $d(w,Z_2^i)$ be the number of neighbors of $w$ in and $Z_2^i\setminus \{v\}$. Also let $f(w)$ be the number of vertices that are connected to $w$ by multiple edges. We consider the following cases:

\textbf{Case a:} $w\in Z_1^i$ then  $\zeta_{i+1}-\zeta_i=-2$.
\\ \textbf{Case b:} $w\in Z_2^i$ and $d(w,Z_2^i)=1$  then $\zeta_{i+1}-\zeta_i=0$.
\\ \textbf{Case c:} $w\in Z_2^i$ and $d(w,Z_2^i)=0$  then $\zeta_{i+1}-\zeta_i=-1$.
\\ \textbf{Case d:} $w\in Z^i \setminus Z_2^i$ then $\zeta_{i+1}-\zeta_i\leq -1+d(w,Z_2^i)+O(f(w))$.

 Summarizing we have, 
\begin{equation}\label{cases}
    \zeta_{i+1}- \zeta_i 
    \begin{cases}
     = -2, & \text{ Case a: probability } (\zeta_i/2m_i)(1+O(m_i^{-1})) .\\
     =  0, & \text{ Case b: probability } p_{2,i}^2  (1+O(m_i^{-1})).\\
=  -1 & \text{ Case c: probability } p_{2,i}(1-p_{2,i})(1+O(m_i^{-1})).\\  
     \leq -1+d(w,Z_2)+O(f(w)) & \text{ Case d: otherwise }
    \end{cases}
  \end{equation}
The contribution of Case d  to $\E(Y_iR_i|\cH_i)$ is at most
\begin{align}\label{problematic}
&\E\big[[-1+d(w,Z_2^i)+O(f(w))]\mathbb{I}(w\in Z^i\setminus Z_2^i) |\cH_i\big]\nonumber
\\ &=-\Pr(w\in Z^i\setminus Z_2^i)+2p_{2,i}\times \frac{3|Z_3|}{2m_i} 
+ 3p_{2,i} \times\frac{4|Z_4|}{2m_i}\nonumber + O\brac{\frac{\lambda_i\log^2m_i}{m_i^{1/2}} + \lambda_i^3}\nonumber\\
&= -\Pr(w\in Z^i\setminus Z_2^i)+p_{2,i}^2\bigg(\lambda_i +\frac{\lambda_i^2}{2}  \bigg)
+ O\brac{\frac{\lambda_i\log^2m_i}{m_i^{1/2}} + \lambda_i^3}.
\end{align}
 For the second line of the above calculation recall that 
$w$ belongs to $Z_j$, $j\geq 3$, with probability $j|Z_j|/2m_i$. Thereafter out of its $j-1$ neighbors (other than $v$), $p_{2,i}$  fraction of them belong to $Z_2$. The terms associated with $Z_j$, $j\geq 5$, have been absorbed by the $O(\lambda_i^3)$ error term.    

To derive the last equality we used \eqref{1}.

Finally observe that \eqref{1} implies
\begin{align}\label{total}
1&= \frac{2|Z_2|+3|Z_3|+4|Z_4|}{2m_i}+\frac{\zeta_i}{2m_i}+ O\brac{\frac{\lambda_i\log^2m_i}{m_i^{1/2}}+\lambda_i^3} \nonumber
\\&= p_{2,i}\bigg(1+ \frac{\lambda_i}{2}+ \frac{\lambda_i^2}{6} \bigg)
+\frac{\zeta_i}{2m_i} + O\brac{\frac{\lambda_i\log^2m_i}{m_i^{1/2}}+\lambda_i^3}.
\end{align}
Therefore \eqref{cases}, \eqref{problematic} and $\Pr(w\in Z^i\setminus Z_2^i)\leq 1-p_{2,i}-\frac{\zeta_i}{2m_i}$ give,
\begin{align}\label{eq:wiri}
\E(Y_iR_i|\cH_i)&\leq  \brac{-\frac{2\zeta_i}{2m_i}- p_{2,i}(1-p_{2,i})} \brac{1+O\bfrac{1}{m_i}}\nonumber\\
& -\bigg(1-p_{2,i}-\frac{\zeta_i}{2m_i}\bigg)+ p_{2,i}^2\bigg(\lambda_i +\frac{\lambda_i^2}{2}  \bigg) + O\brac{\frac{\lambda_i\log^2m_i}{m_i^{1/2}}+\lambda_i^3}\nonumber\\
&= -1-\frac{\zeta_i}{2m_i} 
+ p_{2,i}^2\brac{1+\lambda_i +\frac{\lambda_i^2}{2}  } + O\brac{\frac{\lambda_i\log^2m_i}{m_i^{1/2}}+\lambda_i^3}.\nonumber\\
\noalign{In the last equality we used that $\cB_i$ implies that $\frac{1}{m_i}\ll \frac{\lambda_i\log^2m_i}{m_i^{1/2}}$. 
We now use \eqref{total} to replace -1 by the squared expression to obtain}
& \leq - \bigg[ p_{2,i}\bigg(1+ \frac{\lambda_i}{2}+ \frac{\lambda_i^2}{6} \bigg)
 \bigg]^2   + p_{2,i}^2\bigg( 1+ \lambda_i+  \frac{\lambda_i^2}{2}  \bigg)-  \frac{\zeta_i}{2m_i} + O\brac{\frac{\lambda_i\log^2m_i}{m_i^{1/2}}+\lambda_i^3}\nonumber\\
&\leq - \frac{\lambda_i^2 p_{2,i}^2}{12}+O\brac{\frac{\lambda_i\log^2m_i}{m_i^{1/2}}+ \lambda_i^3}\leq O(m_i^{-1}\log^4 m_i).
\end{align}
For the last inequality we used that in the event $\mathcal{B}_i$ \eqref{1} and \eqref{total} imply that $p_{2,i}=1-o(1)$. In addition, 
\begin{align}
\Pr(R_i=1)&\leq\Pr(\text{Case(a)})+\Pr(\text{Case(c)})+ \Pr(\text{Case(d)})\nonumber\\
& =O\bigg( \frac{\zeta_i}{2m_i}+ p_{2,i}(1-p_{2,i})+\lambda_i\bigg)=O\bigg( \frac{\zeta_i}{2m_i} +\lambda_i\bigg). \label{oneof}
\end{align}
where we have used $1-p_{2,i}=O(\lambda_i)$.

In the event $\cB_i$ we have that $\lambda_i \leq m^{-0.2}$. Hence, if $\zeta_i\leq n^{0.4+\e}$ then $\Pr(R_i=1) \leq m_i^{-0.2}$ and
\begin{align*}
\sum_{j=1}^{m_1} \Pr\brac{\sum_{i=1}^j {R_i}    >n^{0.8+\e/3} }
& \leq \sum_{j=1}^{m_1} \Pr\brac{\sum_{i=1}^j {R_i}\mathbb{I}(m_i>n^{0.8})    >n^{0.8+\e/3}-n^{0.8} } 
\\&\leq m_1 \exp\bigg\{ -\frac{(n^{0.8+\e/3}-n^{0.8}-\sum_{m_i =n^{0.8}}^{m_1} m_i^{-0.2})^2}{2m_1} \bigg\} =o(n^{-10}).
\end{align*}
To obtain the exponential bound, we let $Q_j=\sum_{i=1}^j {R_i}\mathbb{I}(m_i>n^{0.8})$. We have\\
 $\E({Q_j})\leq \sum_{m_i =n^{0.8}}^{m_1} m_i^{-0.2}=O(n^{0.8})$ and then we can use the Chernoff bounds, since our bounds for ${R_i}=1$ hold given the history of the process so far.

It follows that, 
\begin{align}\label{Y}
\sum_{j=1}^{m_1} \Pr&\bigg( \sum_{i=1}^j Y_i \geq n^{0.4+\e/2} \bigg)= \sum_{j=1}^{m_1} \Pr\bigg( \sum_{i=1 }^j Y_i {R_i} \geq n^{0.4+\e/2}\bigg) \nonumber
\\& \leq  \sum_{j=1}^{m_1} \Pr\bigg(\sum_{i=1}^j {R_i}  >n^{0.8+\e/3}\bigg)
+\sum_{j=1}^{m_1} \Pr\bigg( \sum_{i=1 }^j W_i {R_i} \geq n^{0.4+\e/2}\bigg|\sum_{i=1}^j {R_i}  \leq n^{0.8+\e/3} \bigg) \nonumber
\\&\leq o(n^{-8}) + m_1
 \max_{j\leq n^{0.8+\e/3}}  \exp\bigg\{-\frac{\big(n^{0.4+\e/2}-\sum_{m_i =1}^{m_1} m_i^{-1}\log^4 m_i\big)^2}{j\log ^2n}\bigg\} \nonumber
\\& \leq o(n^{-8})+ m_1 \max_{j\leq n^{0.8+\e/3}}  \exp\bigg\{-\frac{\big(n^{0.4+\e/2}-n^{o(1)}\big)^2}{j\log^2 n}\bigg\} =o(n^{-10}).
\end{align}
To obtain the third line we use \eqref{eq:max} which implies that  $|Y_i|\leq \log n$ with probability $1-o(n^{-10})$ for $i\leq \tau$. 
\vspace{3mm}
\\\textbf{Case 3:   $\cC_i$  occurs} \\
At time $\sigma'=\sigma-1$ we have $|Z^{\sigma'}|\lambda_{\sigma'}\geq  \log^4 n$ and hence the estimates \eqref{1} hold.
Thereafter $|Z^{\sigma'}-Z^{\sigma}|,|m_{\sigma'}-m_{\sigma}| = O( \Delta(G_{\sigma'}))$. The maximum degree of $\Delta(G_{\sigma'})$ is bounded  by $\log n$ with probability $1-o(n^{-10})$. 
At time $\sigma$ we have $Z^{\sigma}\lambda_{\sigma}<  \log^4 n$ hence $\lambda_{\sigma'}\leq \frac{2\log^4n}{m_{\sigma'}}$ and so  subsequently for $i\geq \sigma$ we have
\begin{equation}\label{smallsets}
|Z_3^i|=O(\log^4n)\text{ and }Z_{j}^i=\emptyset\text{ for }j\geq 4.
\end{equation}
Given the above we replace \eqref{total} by
\begin{equation}\label{total1}
1=p_{2,i}+\frac{\zeta_i}{2m_i}+O\bfrac{\log^4n}{m_i}.
\end{equation}
Following this just as in \eqref{Y} we get
$$\sum_{j=1}^{m_1} \Pr\bigg( \sum_{i=1}^j Z_i '\geq n^{0.4+\e/2} \bigg)=o(n^{-10}).$$
The above analysis and equation \eqref{goaly} shows that with probability $1-o(n^{-10})$ 
$$\min\{\zeta_i,n^{0.4+\e}\}\leq \log^2 n +4n^{0.4+\e/2} <n^{0.4+0.9\e}.$$
Hence with probability $1-o(n^{-9})$  there does not exist $i<\tau$ such that $\zeta_i> n^{0.4+\e}$. And this therefore completes the proof of Lemma \ref{lem:zeta}.

\end{appendix}
\end{document}